\documentclass[11pt,final]{amsart}
\usepackage{verbatim}
\usepackage[ruled,linesnumbered,lined,boxed,commentsnumbered]{algorithm2e}
\usepackage{latexsym, graphicx, epsfig}
\usepackage{pgfplots}
\usepackage{url}
\usepackage{color}
\usepackage{array}
\usepackage{amsfonts,amsmath,amssymb}
\usepackage{verbatim}
\usepackage{subfigure}
\usepackage{graphicx}

\usepackage[utf8]{inputenc} 
\usepackage[T1]{fontenc}

\usepackage{cite}

\def\mb{\mathbf}

\renewcommand{\epsilon}{\varepsilon}

\SetKwInOut{Offline}{Offline}
\SetKwInOut{Online}{Online}

\usepackage{cite}

\begin{document}

\title{A Bi-fidelity Stochastic Collocation method for 
transport equations with diffusive scaling and multi-dimensional random inputs}

\author{Liu Liu\address{Department of Mathematics, The Chinese University of Hong Kong, Shatin, N.T., Hong Kong SAR }\email{lliu@math.cuhk.edu.hk} 
\and Lorenzo Pareschi\address{Department of Mathematics and Computer Science, University of Ferrara, Via Machiavelli 30, Ferrara,
44121, Italy}\email{lorenzo.pareschi@unife.it}
      \and
             Xueyu Zhu\address{Department of Mathematics, University of Iowa, Iowa City, IA 52242, USA}\email{xueyu-zhu@uiowa.edu}}

\begin{abstract}
In this paper, we consider the development of efficient numerical methods for linear transport equations with random parameters and under the diffusive scaling. 
We extend to the present case the bi-fidelity stochastic collocation method introduced in \cite{NGX14,ZNX14,ZX17}. For the high-fidelity transport model, the asymptotic-preserving scheme \cite{JPT2} is used for each stochastic sample. 
We employ the simple two-velocity Goldstein-Taylor equation as low-fidelity model to accelerate the convergence of the uncertainty quantification process. The choice is motivated by the fact that both models, high fidelity and low fidelity, share the same diffusion limit. Speed-up is achieved by proper selection of the collocation points and relative approximation of the high-fidelity solution. 
Extensive numerical experiments are conducted to show the efficiency and accuracy of the proposed method, even in non diffusive regimes, with empirical error bound estimations as studied in \cite{GZJ20}. 
\end{abstract}
\maketitle
\vskip -.5cm
\noindent \subsubsection*{Keywords.}{\small
 transport equations, uncertainty quantification, stochastic collocation, bi-fidelity method, 
Goldstein-Taylor model, diffusive scaling, asymptotic-preserving schemes}
\vskip +.5cm

\pagestyle{myheadings}
\thispagestyle{plain}
\markboth{L.LIU, L. PARESCHI AND X. ZHU}
     {BI-FIDELITY STOCHASTIC COLLOCATION FOR TRANSPORT EQUATIONS}


\section{Introduction}
\label{sec:Intro}

Transports equations, which describe how kinetic particles make collisions and absorption through a material medium while evolving in time, arises in many applications, such as atmosphere and ocean modeling \cite{coakley2014atmospheric,stamnes2017radiative,zdunkowski2007radiation}, astrophysics \cite{peraiah2002introduction},  nuclear physics \cite{CZ,LK}, biology \cite{Had,Per} and epidemiology \cite{Ber,boscheri2020}. Such problems usually involve a mixed range of length scales characterized by the Knudsen number, defined as the ratio of the mean free path over a typical length scale such as the size of the spatial domain. 
To resolve multiple scales, carefully designed numerical methods are usually needed, see for example \cite{jin2010asymptotic, Acta,boscarino2013implicit,JPT2} for multiscale kinetic problems.  

Besides the multiscale challenges, practical applications of the linear transport model usually contain uncertainties \cite{JH1, JH2, jin2015asymptotic, Ber2}. For example, the scattering cross-section in the collision operator is usually extracted from data, or in some cases we have only a rough estimate of the initial data. Such uncertainties could compromise the predictive capability of the underlying model. Efficient uncertainty quantification (UQ) for such problem becomes critical for improving the reliability of numerical predictions in 
real-life applications, see example \cite{L-Review,JinL-Review,survs} for an overview of broader UQ problems for kinetic and related models.
Among many numerical methods in UQ, stochastic collocation (SC) method has shown its competitiveness in many practical applications due to its intrinsic simplicity and non-intrusive nature. There has been many work in this direction developed in recent years, see for examples \cite{Babuska, Roberts, NXiu, Webster, GWZ14, Xiu07, Schwab} and the references therein. Despite the successful development of SC methods, the high cost of a single high-fidelity run together with the number of required high-fidelity simulation runs, can still render it computationally infeasible for large-scale applications with high-dimensional random parameter space. 

Alternatively, many efforts have been devoted to developing multi-fidelity methods by different communities from various perspectives to mitigate the computational cost \cite{giles2008multilevel,perdikaris2017nonlinear,peherstorfer2016optimal,eldred2017multifidelity,NGX14,ZNX14,ZX17}. See also the recent survey \cite{peherstorfersurvey}.
By combining an inaccurate but cheap low-fidelity solver and an accurate yet expensive solvers in clever ways, multi-fidelity ideas have shown effectiveness in reducing the number of high-fidelity samples. Specifically, within the context of kinetic equations, multi-scale control variate and multilevel variance reduction Monte Carlo methods for kinetic equations with uncertainties has been studied to reduce the sampling cost in \cite{DP19,DP20,HPW20,PTZ}{\color{black}, particularly for problems close to the fluid regime or in the grazing limit}.  A bi-fidelity approximation constructed by combining the compressible Euler system and Boltzmann equations can approximate well the high-fidelity solutions at a much reduced computational cost \cite{liu2020bi,DP19}.

In this work,  we are interested in adapting the bi-fidelity method developed in \cite{NGX14, ZNX14, ZX17} to efficiently approximate the high-fidelity statistics of the multiscale linear transport equation with multi-dimensional random parameters. One key component of the bi-fidelity method is the choice of a proper low-fidelity model, which is crucial for the approximation quality of the bi-fidelity surrogate. In this work, we consider the simple two-velocity Goldstein-Taylor (GT)  model \cite{Gol,Tay} as our low-fidelity model. This choice is motivated by the following observation: after the even and odd parity formulation of the transport equation, the general linear transport equation can be reformulated as a similar form of the discrete-velocity GT model, which can be regarded as a good approximation with discrete velocity variables, and more importantly, 
with significantly less computational cost \cite{JPT2}. Additionally the two models, high-fidelity and low-fidelity, share the same diffusion limit.
We demonstrate that the bi-fidelity approximations based on this low-fidelity model could produce reasonably satisfactory results across a large range of regimes, from the kinetic to diffusive regime. To our knowledge, this is the first manuscript in which a multi-fidelity approach has been applied in the context of diffusive scaling.

The rest of the manuscript is organized as follows. In Section \ref{sec:hi-fi model}-\ref{sec:model}, we introduce the linear transport equation with random inputs and the uncertain discrete-velocity GT model. The former is the high-fidelity model (as it is our goal problem to be solved), whereas an equivalent reformulation of the GT model is served as our low-fidelity model. In Section \ref{sec:pre}, we first give a brief review of the bi-fidelity approach \cite{NGX14, ZNX14} and then recall the corresponding empirical error bound estimations. In Section \ref{sec:NumExample}, extensive numerical experiments are shown to illustrate the effectiveness and efficiency of our proposed method, where kinetic, fluid and mixed regimes are all carefully examined. We end the manuscript with some concluding remarks and future perspectives in Section 6.

\section{Transport equations with random input and diffusive scaling}
\label{sec:hi-fi model}
We consider an one-dimensional linear transport equation with
random scattering coefficients under the diffusive scaling.
Let $f(t,x,v)$ be the probability density distribution of particles at
position $x$, time $t$, and with $v\in(-1,1)$ the cosine of the angle
between the particle 
velocity and its position variable. Then $f$ is governed by the following
linear transport equation \cite{CZ}:
\begin{equation}
\label{pde_transport1d}
\epsilon \partial_t  f + v \partial_x f = \frac{\sigma(x,z)}{\epsilon}
\left[\frac{1}{2}\int_{-1}^1 f(v')\, dv' -f\right], 
\end{equation}
where $\epsilon$ is the Knudsen number and  $\sigma(x,z)$ is the random scattering coefficient and $z \in I_z$ with $I_z$ the domain of the random parameter space. In this project,  
we ignore the absorption and source terms, which can also introduce uncertainties. The treatment of these terms does not add numerical difficulties thereby will be neglected here for the ease of presentation.

To understand its diffusion limit, we first split equation \eqref{pde_transport1d} into
two equations for $v>0$: 
\begin{equation}
\label{pde_transport1d_sp}
\begin{split}
&\epsilon \partial_t  f(v) + v \partial_x f(v) =
\frac{\sigma(x,z)}{\epsilon}\left[\frac{1}{2}\int_{-1}^1f(v')\,
  dv'-f(v)\right], \\[4pt]
&\epsilon \partial_t f(-v) - v \partial_x f(-v) =
\frac{\sigma(x,z)}{\epsilon}\left[\frac{1}{2}\int_{-1}^1f(v')\,
  dv-f(-v)\right]. 
\end{split}
\end{equation}
By the even-odd decomposition \cite{JPT2}, we introduce the even and odd parities 
\begin{equation}
\begin{split}
r(t,x,v) &= \frac{1}{2}[f(t,x,v) + f(t,x,-v)], \\[4pt]
j(t,x,v) &= \frac{1}{2\epsilon}[f(t,x,v) - f(t,x,-v)].
\end{split}
\end{equation}
The system \eqref{pde_transport1d_sp} can then be reformulated by the following system: 
\begin{equation}
\label{pde_transport1d_rj}
\left\{
\begin{split}
&\partial_t r + v \partial_x j = \frac{\sigma(x,z)}{\epsilon^2}(\overline{r}-r), \\[4pt]
&\partial_t j + \frac{v}{\epsilon^2} \partial_x r =
-\frac{\sigma(x,z)}{\epsilon^2}j, 
\end{split}
\right.
\end{equation}
where 
\begin{equation}
\label{r_integral}
\overline{r}(t,x)=\int_0^1 rdv.
\end{equation}
As $\epsilon\rightarrow 0^+$, \eqref{pde_transport1d_rj} yields
\[r=\overline{r}, \qquad j=-\frac{v}{\sigma(x,z)}\partial_x \overline{r}.\]
Substituting it into system (\ref{pde_transport1d_rj}) and integrating over $v$, one gets the limiting diffusion equation (\cite{LK, BSS}) as follows: 
\begin{equation}
\label{transport1d_diff}
\left\{
\begin{split}
&j=-\frac{v}{\sigma(x,z)}\partial_x \overline{r}, \\[4pt]
&\partial_t \overline{r} =\partial_x \left[\frac{1}{3\sigma(x,z)} \partial_x\overline{r}\right].
\end{split}
\right.
\end{equation} 
Solving the model \eqref{pde_transport1d_sp} brings many numerical challenges due to the complexity of velocity integral operators on the right-hand-side, the stiffness of the source terms, and in particular for the case with high-dimensional random parameters \cite{JH2, JPT2}. In this regard, we desire to adapt the bi-fidelity approximation method developed in \cite{NGX14, ZNX14} to the present case in order to mitigate the computational cost of a standard Monte Carlo sampling approach.

\section{The Low-fidelity Goldstain-Taylor model}
\label{sec:model}
In this section, we introduce a simple discrete-velocity model
subject to random inputs and discuss the motivation of choosing it as our low-fidelity model for the linear transport equation within the framework of bi-fidelity approximation.

The one-dimensional Goldstein-Taylor (GT) model \cite{Gol,Tay} is given by
\begin{equation}
\label{GT}
\left\{
\begin{split}
& \partial_t u + \frac{1}{\epsilon}\partial_x u = \frac{\sigma(x,z)}{2\epsilon^2}(v-u), \\[4pt]
& \partial_t v - \frac{1}{\epsilon}\partial_x v = \frac{\sigma(x,z)}{2\epsilon^2}(u-v). 
\end{split}
\right.
\end{equation}
Here we assume a random scattering coefficient $\sigma(x,z)$, that already in itself characterizes a model with many interesting applications \cite{GottliebX_CICP08}. 
We introduce the macroscopic variables: the mass density $\rho$ and the flux $s$,
$$ \rho = u + v, \qquad s = \frac{u-v}{\epsilon}, $$
then the GT model \eqref{GT} is equivalent to the following system: 
\begin{equation}
\label{pde_transport1d_uv}
\left\{
\begin{split}
&\partial_t \rho +  \partial_x s = 0, \\[4pt]
&\partial_t s + \frac{1}{\epsilon^2}\partial_x\rho =
-\frac{\sigma(x,z)}{\epsilon^2}s. 
\end{split}
\right.
\end{equation}
The above system is the analogous of \eqref{pde_transport1d_rj} for the high-fidelity model and in the diffusion limit $\epsilon\rightarrow 0$, system \eqref{pde_transport1d_uv}
can be approximated by the heat equation to the leading order, with random diffusion coefficient $\sigma(x,z)$: 
\begin{equation}
\label{pde_transport1d_sp_lowfi2_diff}
\left\{
\begin{split}
&s = -\frac{1}{\sigma(x,z)}\partial_x \rho, \\[4pt]
&\partial_t \rho = \partial_{x} 
\left[\frac{1}{\sigma(x,z)} \partial_x \rho \right].
\end{split}
\right.
\end{equation}

Comparing \eqref{transport1d_diff} and \eqref{pde_transport1d_sp_lowfi2_diff}, both systems look similar except for the magnitude of the diffusion coefficient on the right hand-side. 
{If one sets the diffusion coefficient in \eqref{pde_transport1d_sp_lowfi2_diff} be $\sigma_{\text{GT}}$ and that in \eqref{transport1d_diff} be $\sigma_{\text{LTE}}$, then by assuming $\sigma_{\text{GT}}=\frac{1}{3}\sigma_{\text{LTE}}$ the two models share the same diffusion limit.}
It is worth noting that 
the GT model \eqref{GT} can be regarded as a discrete-velocity kinetic model of the linear transport equation: $u$ defines the density of particles traveling with velocity $1$, 
whereas $v$ that of particles traveling in the reverse direction with velocity $-1$. Besides, the GT model has much cheaper computational cost, yet shares the same limiting diffusion equations with the linear transport model as $\epsilon\to 0$. Motivated by these observations, we employ the equivalent formulation of the GT equation, that is, system \eqref{pde_transport1d_uv} as our low-fidelity model. We refer to \cite{Lions1997} for rigorous results concerning the diffusion limit of two-velocity models and extensions to nonlinear diffusion coefficients.

\section{The bi-fidelity stochastic collocation approach} 
\label{sec:pre}

In this section, we briefly review the bi-fidelity method. Assume the expensive high-fidelity model and the cheap low-fidelity model are available to generate the high-fidelity solution $u^H(z)$ and the low-fidelity solution $u^L(z)$ respectively, for any given parameter $z$. The main idea of the bi-fidelity approximation in \cite{NGX14, ZNX14} is to approximate the high-fidelity solution $u^H(z)$ by the following expansion: 
\begin{equation}
 u^B(z) = \sum_{k=1}^nc_k(z) u^H(z_{k}),
\end{equation}
where $n$ is the  number of the selected high-fidelity solutions in the parameter space. If $n$ is small and the coefficient $c_k(z)$ can be efficiently and accurately approximated, an efficient bi-fidelity approximation can be constructed. To achieve this goal, two major questions need to be addressed:
(1) how to select the collocation points $z_k$ so that the total number of high-fidelity samples $n$ is small? (2) how to approximate the expansion coefficients $c_k(z)$ properly for any given parameter $z$? 
\bigskip

{\bf Point Selection}.
It would be computationally infeasible to select the important points among high-fidelity samples due to the cost of the high-fidelity solver. To migrate this cost, we search the important points in the parameter space guided by the low-fidelity model, which is informative yet cheap to evaluate over a large number of points in the parameter space.

Specifically, we denote the candidate set $\Gamma_N=\{z_1, z_2,\hdots, z_N\}$, which is assumed to be large enough to cover the parameter space $Z$.  We shall identify important points iteratively by a greedy approach \cite{NGX14, ZNX14}. Initially, denote $\gamma_0=\{\}$ and assume we have the first $k$ important points  $\gamma_k=\{z_{i_1}, z_{i_2}, \hdots, z_{i_k}\}$ available at the $k$-th iteration. Denote the snapshot matrix $u^L(\gamma_k)= \{u(z)| z \in \gamma_k\}$ and the corresponding spanned solution space $U^L(\gamma_k) = \text{Span}\{u^L(z)|z\in\gamma_k\}$. Then we  pick the $z$ point (from the candidate set $\Gamma_N$) so that  the  corresponding low-fidelity solution is farthest away from the existing spanned solution space $U^L(\gamma_k)$, to be the next sampling point: 
\begin{equation}\label{greedy}
z_{i_{k+1}} = \arg \max_{z\in \Gamma_N} d^L(u^L(z), U^L(\gamma_k)), \quad \gamma_{k+1} = \gamma_k \cup z_{i_{k+1}}, 
\end{equation}
where $d^L(v,W)$ is the distance between a function $v\in u^L(\Gamma_N)$ and the space $W\in u^L(\gamma_k)$. We then continue this process until all $n$ important points $\gamma_n$ are selected. The whole procedure can be efficiently implemented by performing the pivoted Cheloskey decompostion on $u^L(\Gamma_N)$. We refer the reader to more details in \cite{NGX14, ZNX14}. With the selected parameter points $\{\gamma_n\}$, the corresponding high-fidelity approximation space can be constructed by $U^H(\gamma_n) = \text{Span}\{u^H(z)|z\in\gamma_n\}$.
\bigskip

{\bf Approximation of bi-fidelity sample}.
Ideally, one can  compute the expansion coefficients by projecting the high-fidelity data onto the high-fidelity approximation space $U^H(\gamma_n)$, which requires the high-fidelity simulation during the online stage. To migrate this cost, the bi-fidelity approach developed in \cite{NGX14, ZNX14} uses the low-fidelity coefficients $c_k^L(z)$ as an approximation for the high-fidelity coefficients $c_k(z)$. In other words, 
for any given $z$, we shall compute the low-fidelity solution $u^L(z)$ and its low-fidelity coefficients by projecting onto the low-fidelity approximation space $U^L(\gamma_n)$:
\begin{equation}
u^L(z) \approx \mathcal{P}_{U^L(\gamma_n)}u^L(z) = \sum_{k=1}^n c_k^L(z)u^L(z_{k}), \quad z_{k}\in \gamma_n. 
\end{equation}
where the low-fidelity projection coefficients can be computed as follows:
\begin{equation}
{\bf G}^L {\bf c}^L = {\bf f}, \qquad {\bf f}^L = (f_k^L)_{1\leq k\leq n}, \qquad f_k^L = 
\langle u^L(z), u^L(z_{k})\rangle,
 \end{equation}
and
 ${\bf G}^L$ is the Gramian matrix of $u^L(\gamma_n)$, 
\begin{equation}
 ({\bf G}^L)_{ij} =  \left\langle u^L(z_i), u^L(z_j) \right\rangle^L, \qquad 1 \leq i,\, j \leq n, 
\end{equation}
where $\langle\cdot,\cdot\rangle^L$ is the standard inner product associated with $U^L(\gamma_n)$. 

Once $c^L_k(z)$ are computed, it can serve as surrogates for the high-fidelity coefficients shown below in \eqref{UB}. Consequently, the bi-fidelity approximation of the high-fidelity approximation
solution $u^H(z)$ can be constructed by the following: 
\begin{equation}
\label{UB}
 u^B(z) = \sum_{k=1}^nc_k^L(z) u^H(z_{k}).
\end{equation}
To make things clearer, the bi-fidelity approximation of the high-fidelity sample for a given $z$ is summarized in 
{Algorithm ~\ref{BiFi-pod}}.

\begin{algorithm*}[htb]
\caption{Bi-fidelity approximation for a high-fidelity solution at given $z$}
\label{BiFi-pod}
\Offline

Select a sample set $\Gamma_N = \{z_1, z_2, \hdots, z_N\}\subset I_z $.

Run the low-fidelity model $u^L(z_j)$ for each $z_j \in \Gamma_N$.

Select $n$ ``important'' points from $\Gamma_N$. Denote  $\gamma_N=\{z_{i_1}, \cdots z_{i_n} \} \subset\Gamma_N$ and the low-fidelity approximation space by $U^L(\gamma_n)$.

Run high-fidelity simulation at each point in the selected sample set $\gamma_n$. 

\Online

For any given $z$, compute the low-fidelity solution $u^L(z)$ and the corresponding low-fidelity coefficients by projection:
\begin{equation}\label{C-N}
u^L(z) \approx \mathcal{P}_{U^L(\gamma_n)}u^L(z) = \sum_{k=1}^n c_k^L(z)u^L(z_k), \quad z_k \in \gamma_n. \end{equation}
During the online stage, the operator $\mathcal P_{U^L(\gamma_n)}$ in \eqref{C-N} is a projection onto the space $U^L(\gamma_n)$ with projection coefficients $\{ c_k^L\}$ computed by the Galerkin approach: 
\begin{equation}\label{Gc}
{\bf G}^L {\bf c}^L = {\bf f}, \qquad {\bf f}^L = (f_k^L)_{1\leq k\leq N}, \qquad f_k^L = 
\langle u^L(z), u^L(z_{k})\rangle,
 \end{equation}
where ${\bf G}^L$ is the Gramian matrix of $u^L(\gamma_n)$, 
\begin{equation}\label{GM} 
 ({\bf G}^L)_{ij} =  \left\langle u^L(z_i), u^L(z_j) \right\rangle^L, \qquad 1 \leq i,\, j \leq n, 
\end{equation}
with $\langle\cdot,\cdot\rangle^L$ the inner product associated with $U^L(\gamma_n)$. 

Construct the bi-fidelity approximation by applying the same approximation rule as in low-fidelity model:
\begin{equation} u^B(z)  = \sum_{k=1}^n c_k^L(z)u^H(z_k). \end{equation}
\end{algorithm*}

\medskip
{\bf{Bi-fidelity Mean}}. Once we have the bi-fidelity surrogate $u^B(z)$, it is straightforward to employ the Monte Carlo or other quadrature-based methods to compute the statistical moments, such as expectation: 
\begin{equation}\label{mc_moments}
\mathbb{E}[u^H]  \approx \mathbb{E}_M[u^H] =  \sum_{i=1}^{M} w_i u^B(z_i), 
\end{equation}
where $(z_i, w_i)$ are the points and associated weights in the Monte Carlo method or quadrature rules. 
However, this might still require many bi-fidelity solution reconstructions. In \cite{zhu2017multi} an extension to this approach is developed to approximate high-fidelity solution expectations more efficiently. We first compute the low-fidelity sample mean (via the Monte Carlo or quadrature rules):
\begin{equation}
\mu^L = \sum_{i=1}^{M} w_i u^L(z_i),
\end{equation}
then project it on the low-fidelity approximation space $U^L(\gamma_n)$, 
\begin{equation}
\mu^L \approx \mathcal{P}_{U^L(\gamma_n)} \mu^L =\sum_{i=1}^{n} c_i^L u^L(z_i),
\end{equation}
where the expansion coefficients $c^L$ are computed by solving the following linear system:
\begin{equation}
 {\bf G}^L c^L = \mb{g}^L, \quad g^L = \left\langle \mu^L, u^L(z_{j})  \right\rangle^L , \qquad 1 \leq j \leq n, \, z \in \gamma_n.
\end{equation}
with this coefficient $c^L_{k}(z)$, the bi-fidelity approximation of the high-fidelity mean can be constructed as follows:
\begin{equation}
\mu^B = \sum_{k=1}^{n} c_{k}^L u^H(z_{k}),\quad z_k\in \gamma_n. 
\end{equation}
We refer readers to \cite{zhu2017multi} for additional details.

\subsection{An empirical error bound estimation}
\label{sec:errorEstimate}

In practical applications, \emph{a priori} assessment of the model quality and prediction errors is important. 
A previous study \cite{GZJ20} introduced a novel empirical error bound estimation approach with ease of implementation to evaluate the performance of the bi-fidelity surrogates a priori. In this section, we will briefly describe the methodology. 

In \cite{GZJ20}, an important quantity that characterizes the {\it similarity} between the LF and HF models is introduced: 
\begin{equation}
\label{eqn:relativeDist}
R_{s}(\mathbf{z}) = \frac{d^H(\mathbf{v}^H(\mathbf{z}),U^H(\gamma^{k}))}{||\mathbf{v}^H(\mathbf{z})||}\big /\frac{d^L(\mathbf{v}^L(\mathbf{z}),U^L(\gamma^{k}))}{||\mathbf{v}^L(\mathbf{z})||}. 
\end{equation}
For instance, $R_{s}\approx 1$ implies that the LF model is informative enough in the BF reconstruction. With $R_s$ and the observation in \cite[Theorem 1]{GZJ20}, for any given new point $\mb{z}_*$, one has
\begin{equation}
\begin{split}
\label{eq:eb_a}
\frac{||\mathbf{v}^H(\mathbf{z}_*) - \mathbf{v}^B(\mathbf{z}_*)||}{||\mathbf{v}^H(\mathbf{z}_*)||}
&\leq \frac{d^L(\mathbf{v}^L(\mathbf{z}_*),U^L(\gamma^{k}))}{||\mathbf{v}^L(\mathbf{z}_*)||} R_s(\mathbf{z}_*)\\
&\times\Big(1+\frac{|| P_{U^H({\gamma^k})}\mathbf{v}^H(\mathbf{z}_*) -  \mathbf{v}^B(\mathbf{z}_*)||}{ d^H(\mathbf{v}^H(\mathbf{z}_*),U^H(\gamma^{k}))}\Big).
\end{split}
\end{equation}

To remove the dependency of the new HF sample  $\mathbf{v}^H(\mathbf{z}_*)$ on the above right-hand side, one uses $\mathbf{z}_{k+1} \in \gamma_{k+1}$ as the testing points served as an error surrogate for the BF approximation in the entire parameter space. If the LF and HF models are similar (i.e., $R_s\approx 1$), one can choose some proper constants $c_1$ and $c_2$, such that for the first $k+1$ pre-selected important points $\gamma_{k+1}$, 
\begin{equation}
\label{eq:eb1}
\begin{split}
    \frac{||\mathbf{v}^H(\mathbf{z}_*) - \mathbf{v}^B(\mathbf{z}_*)||}{||\mathbf{v}^H(\mathbf{z}_*)||}
&\leq
\frac{d^L(\mathbf{v}^L(\mathbf{z}_*),U^L(\gamma^{k}))}{||\mathbf{v}^L(\mathbf{z_*})||}\\
&\times \Big[c_1+c_2\frac{|| P_{U^H({\gamma^k})}\mathbf{v}^H(\mathbf{z}_{k+1}) -  \mathbf{v}^B(\mathbf{z}_{k+1})||}{ d^H(\mathbf{v}^H(\mathbf{z}_{k+1}),U^H(\gamma^{k}))}\Big].
\end{split}
\end{equation}
It can be seen that one only needs the LF data and the first $k+1$ pre-selected HF samples. 

Another important quantity introduced in \cite{GZJ20}, called $R_e$, indicates the approximation quality of the BF approach, 
\begin{equation}
\label{eqn:rError}
R_{e}(\mathbf{z}) := \frac{|| P_{U^H({\gamma^k})}\mathbf{v}^H(\mathbf{z}) -  \mathbf{v}^B(\mathbf{z})||}{d^H(\mathbf{z},U^H(\gamma^{k}))}, 
\end{equation}
which describes the balance between the in-plane error and the relative distance. 
When $R_{e}$ is large, one should stop collecting new HF samples. 
Take the expectation of \eqref{eq:eb1} on both sides with respect to $\mb{z}_*$, we get
\begin{align}
\label{eq:eb2}
\begin{split}
    \mathbb{E}\Big[\frac{||\mathbf{v}^H(\mathbf{z}_*) - \mathbf{v}^B(\mathbf{z}_*)||}{||\mathbf{v}^H(\mathbf{z}_*)||}\Big]
& \leq \mathbb{E}\Big[\frac{d^L(\mathbf{v}^L(\mathbf{z}_*),U^L(\gamma^{k}))}{||\mathbf{v}^L(\mathbf{z}_*)||}\left(c_1+c_2R_{e}(\mathbf{z}_{k+1})\right)\Big]. \\
& \leq \Big[\max_{\mb{z}_*}\frac{d^L(\mathbf{v}^L(\mathbf{z}_*),U^L(\gamma^{k}))}{||\mathbf{v}^L(\mathbf{z}_*)||}\Big]
\left(c_1+c_2R_{e}(\mathbf{z}_{k+1})\right)
\end{split}
\end{align}

We acknowledge that the above error bound estimation, though not really rigorous, is a useful quantity to access the quality of the BF approximation in practical applications. In the following numerical experiments, {our empirical results suggest that this error bound estimation is effective if the constants $c_1$ and $c_2$ are set to be  1}.

\section{Numerical Examples}
\label{sec:NumExample}
In this section, we present several numerical examples to illustrate
the effectiveness of our method. 
To examine the accuracy, two metrics are used to quantify the errors: the
differences in the mean solutions and in the corresponding standard 
deviation compared with the reference solutions, with $L^2$ norm in space, 
\begin{equation}
 e_{mean}(u^B) = \left\|\mu^B - \mathbb{E}[u^H] \right\|_{L^2},  \qquad
 e_{std}(u^B) = \left\|\mathbb{\sigma}[u^B] - \mathbb{\sigma}[u^H] \right\|_{L^2},
\end{equation}
where $\mu^B, \sigma^B$ are the bi-fidelity approximation for the high-fidelity mean and standard deviation, respectively. Here $\mathbb{E}[u^H],\sigma^H$ are the corresponding high-fidelity sample mean and sample standard deviation that served as reference solutions. 

Without loss of generality, the $d$-dimensional random variable $\mathbf{z}=\{z_1, \cdots, z_d\}$ is assumed to follow the uniform distribution on $[-1,1]^d$ in all our numerical tests, and the dimension of the random parameter is chosen as $d=5$ for simplicity. Higher dimensional random spaces can be treated similarly, note however that already the $d=5$ case would make the corresponding stochastic Galerkin approach \cite{JH2} extremely expensive from a computational viewpoint.

To compute the reference solutions for the mean and standard deviation of the high-fidelity quantities of interests, we use the high-order stochastic collocation method over 5-dimensional 
sparse quadrature points with $5$-level Clenshaw-Curtis rules, i.e., evaluated on $2243$ quadrature points (cf.,\cite{XiuH_SISC05}).  

For the high-fidelity (HF) solver, at each given sample we employ the AP scheme \cite{JPT2} developed for the deterministic linear transport equation \eqref{pde_transport1d} under the diffusive scaling. 
The standard 16-points Gauss-Legendre quadrature set is used for the velocity space to compute $\bar{r}$ in \eqref{r_integral}.
For the low-fidelity (LF) solver, we use the deterministic AP method \cite{JPT1} to solve the linear Goldstein-Taylor model \eqref{pde_transport1d_uv}. Similar results are obtained using IMEX Runge-Kutta approaches of the type proposed in \cite{boscarino2013implicit} for both the high fidelity and the low fidelity models. 

The spatial and temporal discretizations are specified in each of the following tests. Periodic boundary conditions are considered, if not specified. 
Regarding the comparisons on CPU time, since for different numerical examples, different temporal and spatial discretizations in HF and LF models are used, we test and compare the CPU computational time roughly by estimating the HF costs about 20 times more than the LF solvers.

\subsection*{Test 1: Uncertain cross-section}
We first consider the linear transport equation (\ref{pde_transport1d_rj}) with the random cross-section coefficient given by
\begin{equation}
\label{Sigma}
\sigma (x,z) = 1 + \sigma \sum_{i=1}^d \frac{1}{(i\pi)^2}\cos{(2\pi i x)} z_i,
\end{equation}
where $\sigma = 4$ and $d=5$. 
The initial conditions are 
\begin{equation}
r(x,v,z,0) =0, \qquad j(x,v,z,0)=0,
\end{equation}
and the boundary conditions are (see \cite{JPT2})
\begin{equation}
\sigma j=-vr_x,
\end{equation}
thus one gets {\color{black}for $v>0$}, 
\begin{equation}
\left. r-\frac{\epsilon}{\sigma} v r_x \right|_{x=0} =1, \qquad
\left. r+
\frac{\epsilon}{\sigma} v r_x \right|_{x=1} = 0. 
\end{equation}

\begin{figure}[tb]
\centering
\includegraphics[scale=0.33]{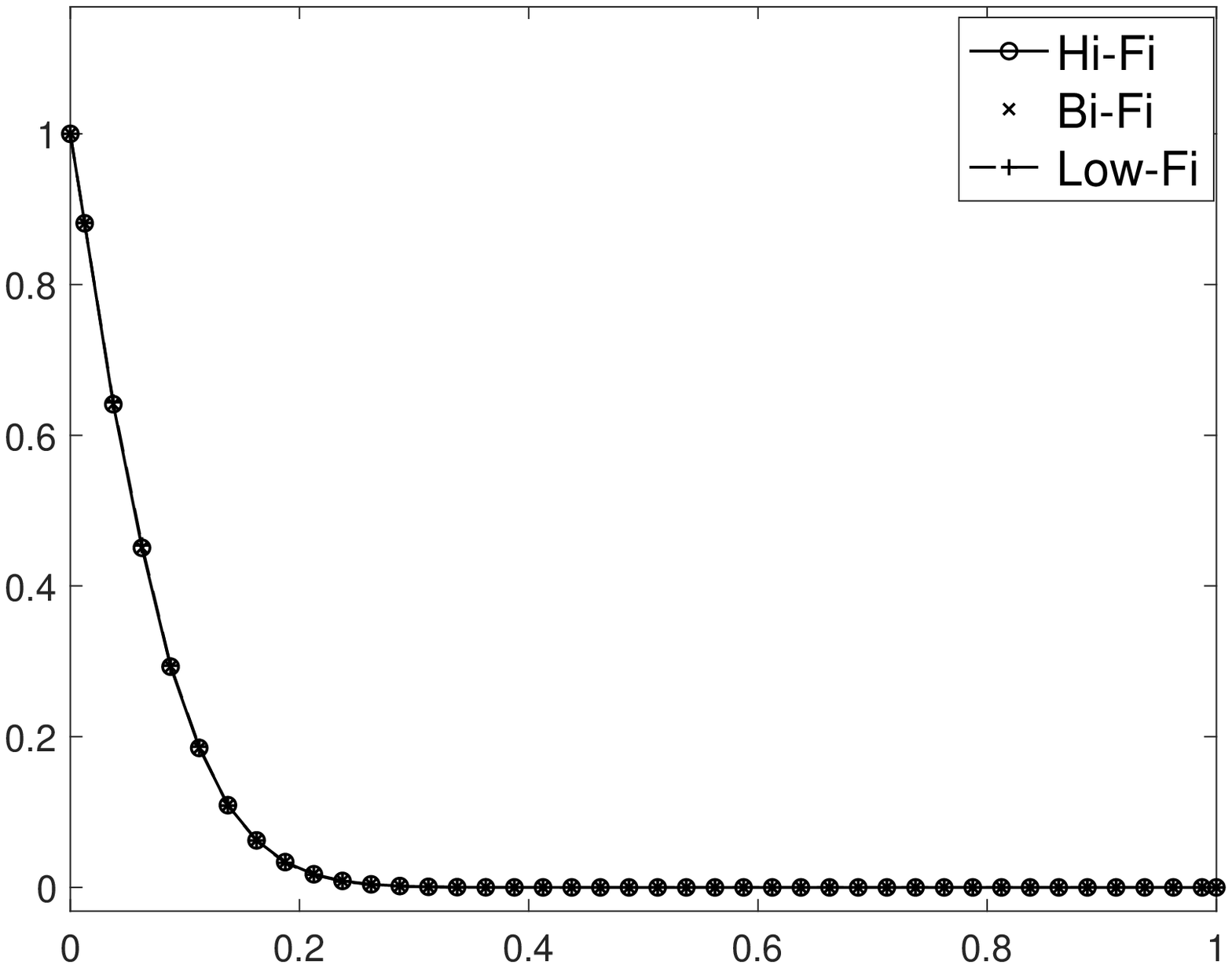}
\includegraphics[scale=0.33]{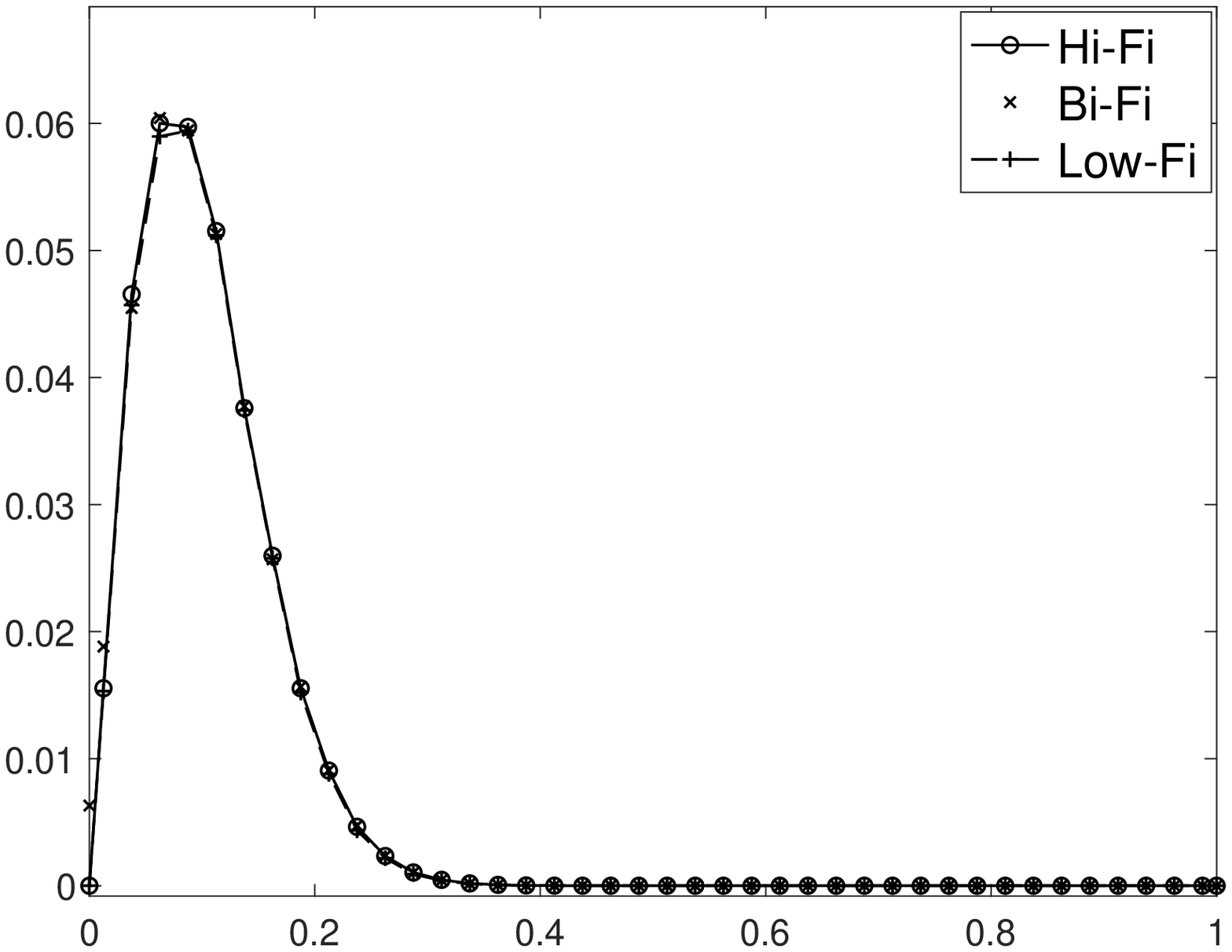}
\includegraphics[scale=0.33]{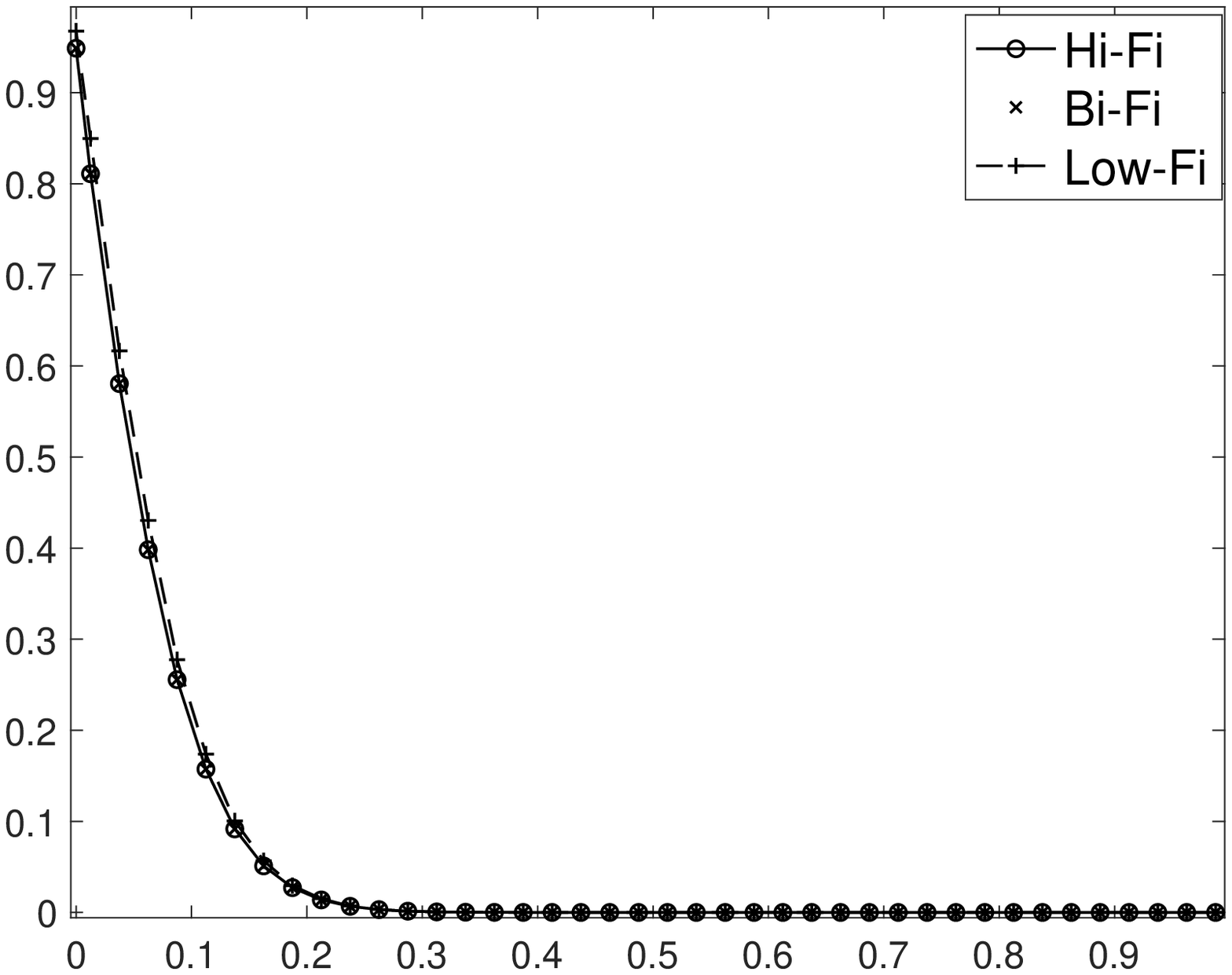}
\includegraphics[scale=0.33]{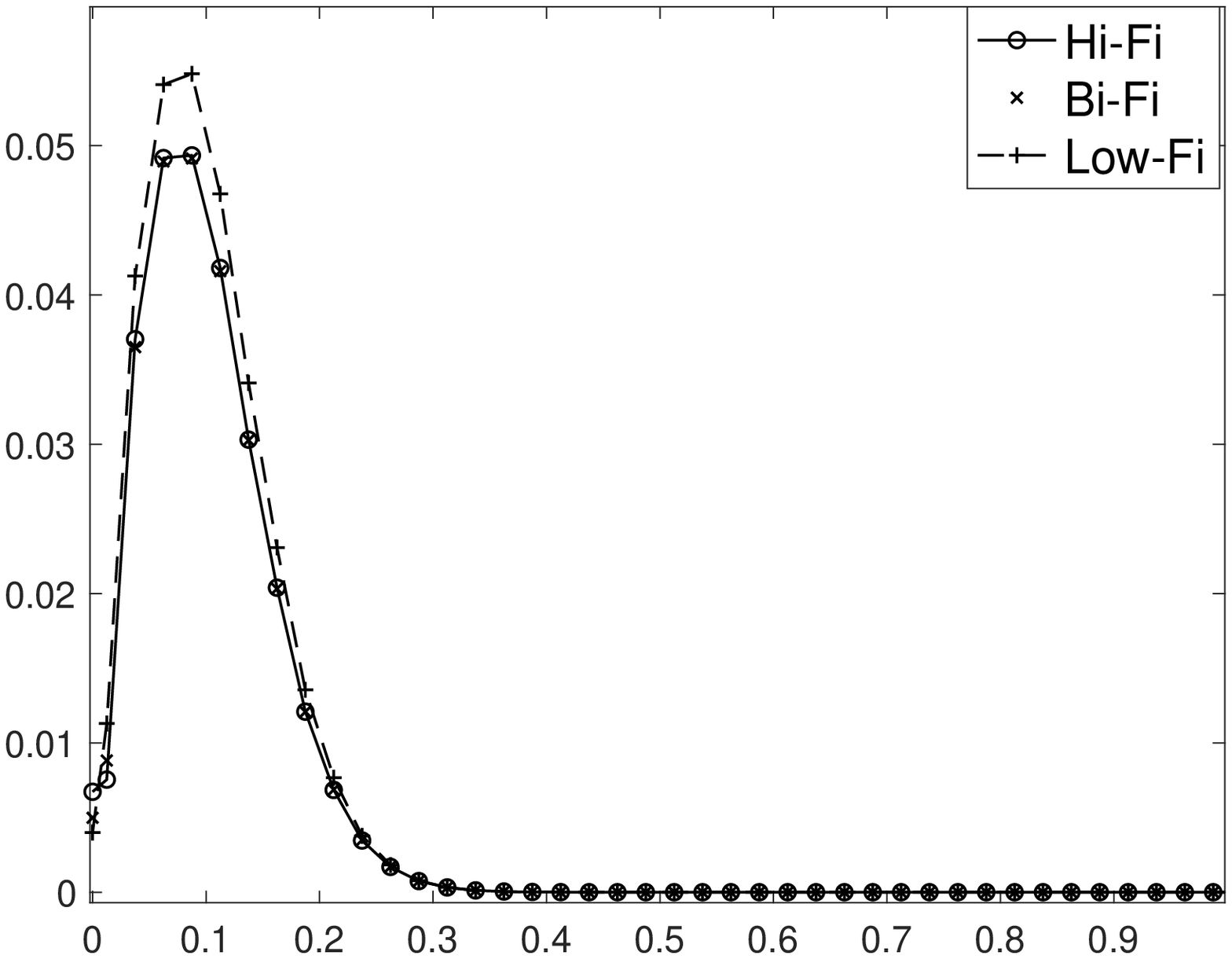}
\caption{Test 1. The mean (left) and standard deviation (right)
  of $\overline{r}$ at $\epsilon=10^{-8}$ (first row), and $\epsilon=10^{-2}$ (second row), obtained by $12$ high-fidelity runs and the sparse grid method with $2243$ quadrature points (crosses).}
\label{1dtransport_comparision_multiD_a}
\end{figure}

\begin{figure}[htb]
\centering
\includegraphics[scale=0.3]{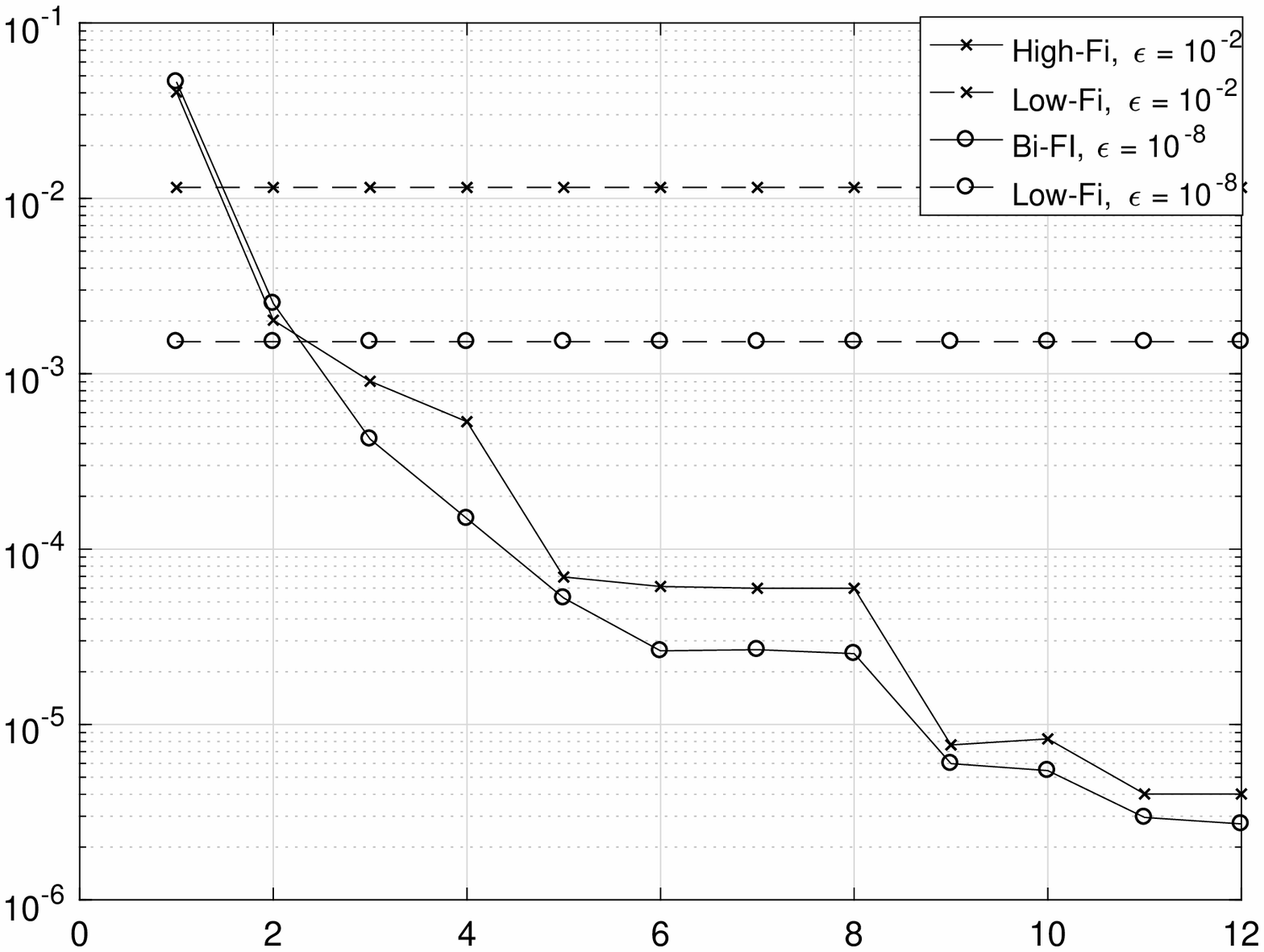}
\includegraphics[scale=0.3]{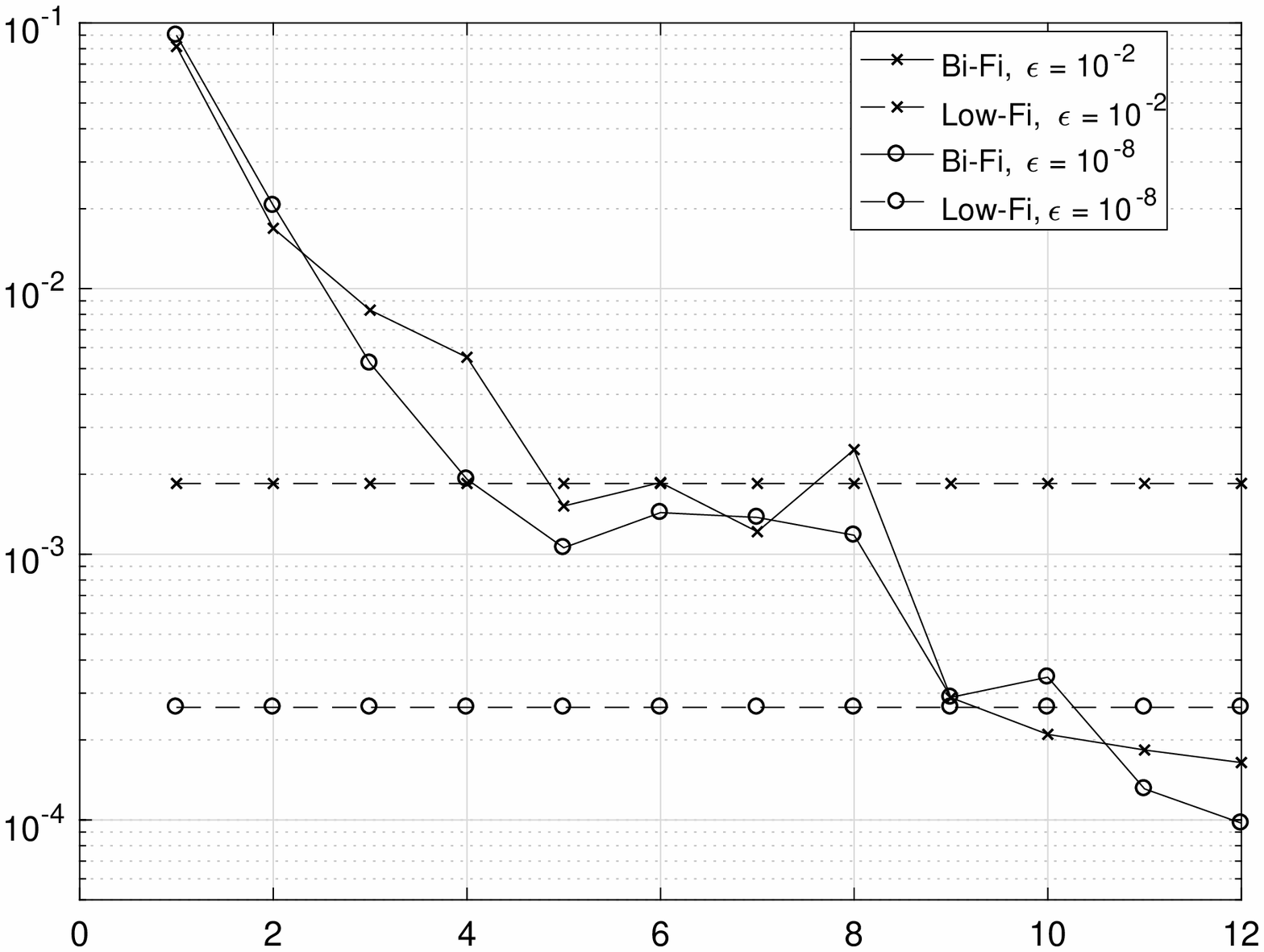}
\caption{Test 1. Error for the bi-fidelity approximation of mean (left) and standard deviation (right) of $\overline{r}$ at $\epsilon=10^{-8}$ and $\epsilon=10^{-2}$, with respect to $r$ high-fidelity runs. }
\label{1dtransport_convergence_multiD_b} 
\end{figure}

\begin{figure}[htb]
\centering
\includegraphics[scale=0.33]{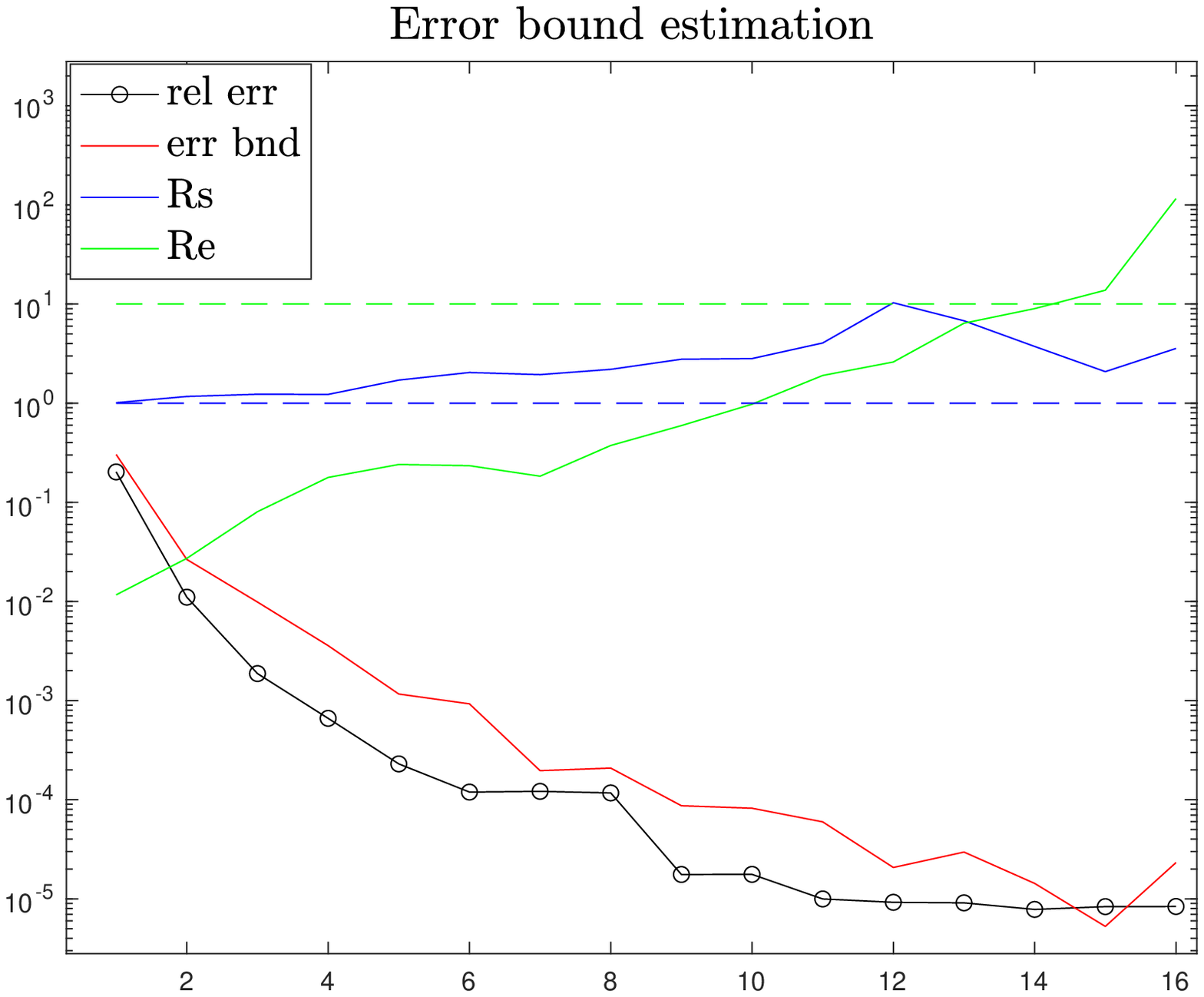}
\includegraphics[scale=0.34]{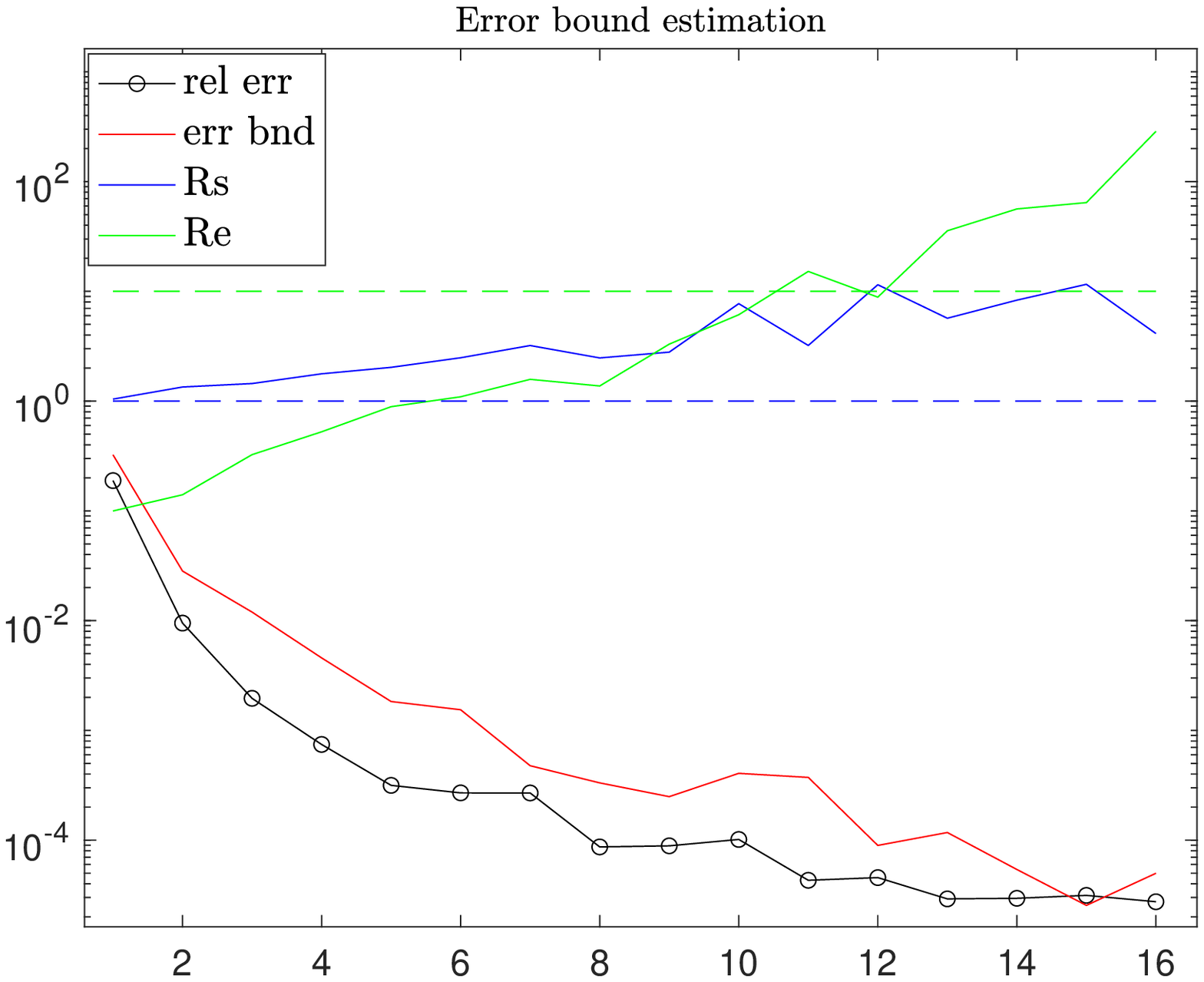}
\caption{Error bound estimation for Test 1. 
Left: $\epsilon=10^{-8}$. Right: $\epsilon=10^{-2}$.}
\label{Test1_Error}
\end{figure}

In this example, we set $\Delta x=0.025$, $\Delta t=2/3\times 10^{-4}$ in the HF solver, while $\Delta x=0.025$, $\Delta t=2\times 10^{-4}$ in the LF solver. 
The mean and standard deviation of the BF solution $\overline{r}$ at time $T=0.01$ for $\epsilon=10^{-8}$ and $10^{-2}$ are shown in Fig.~\ref{1dtransport_comparision_multiD_a}. One can see that the mean and standard deviation of $\overline{r}$ obtained by BF solutions (based on 12 high-fidelity samples) match perfectly well with the HF solutions. In contrast, the LF solutions are slightly off, especially near the left boundary. 

In Fig.~\ref{1dtransport_convergence_multiD_b}, we show
the errors of the mean (circle) and standard deviation (stars) of the BF approximation of
$\overline{r}$ with respect to the number of HF runs for different $\epsilon$. While the baseline
 low-fidelity approximations have the error larger than $\mathcal{O}(10^{-3})$,  the errors can be reduced to $\mathcal{O}(10^{-5})$ for the mean and $\mathcal{O}(10^{-4})$ for the standard deviation of HF solutions with $12$ selected HF samples used in the BF approach. {\color{black}In contrast, the high-fidelity reference solution used $2243$ quadrature points in the sparse grid method, the saving of computational cost is quite noticeable in our BF method.}
  
Besides, one interesting observation is that as $\epsilon$ increase, even though the LF model deviates from the HF model, the bi-fidelity solutions seem to
show similar convergence trend in both cases, indicating a uniform (in $\epsilon$) accuracy properties of the bi-fidelity approximation. Some theoretical studies for a general class of multiscale kinetic models was conducted on this topic, see \cite{Bi-Error}.

From Fig.~\ref{Test1_Error}, one can see that 1) the model similarity $R_s$ is close to $1$ and $R_e$ is less than $10$  when the number of HF samples  is small; 
2) the empirical error bound estimations (red solid lines) is able to bound well the true errors, and adding more HF training samples can further reduce the BF errors; 3) until when the number of HF simulations reaches $r=12$, the similarity metric $R_s$ becomes more fluctuating and $R_e$ also increase significantly as $r$ further grows, which indicates that collecting additional HF samples does not help further to improve the approximation quality of our BF surrogate.

\subsection*{Test 2: Uncertain cross-section and initial data}
In this test, we consider the random cross-section coefficient given by \eqref{Sigma} as well as uncertain initial data: 
\begin{equation}
\label{Test2_IC}
f(t=0,x,v,z)=\rho_0\exp\left(-\left(\frac{v-0.5}{T_0}\right)^2\right)+ 
\rho_1 \exp\left(-\left(\frac{v+0.75}{T_1}\right)^2\right), 
\end{equation}
where 
\begin{align*}
& \rho_0(x,z) = 1 + 3 \sum_{k=1}^d \frac{1}{(k\pi)^2} \sin(2\pi k x)z_k, 
\qquad T_0(x,z)=\frac{5+2\cos(2\pi x)}{20}\left(1+0.6 z_1 \right), \\[2pt]
& \rho_1(x,z) = 1 + 2 \sum_{k=1}^d \frac{1}{(k\pi)^2} \cos(2\pi k x)z_k, 
\qquad T_1(x,z) = 0.5+0.2\cos(2\pi x)z_2. 
\end{align*}

\begin{figure}[htb]
\centering
\includegraphics[scale=0.32]{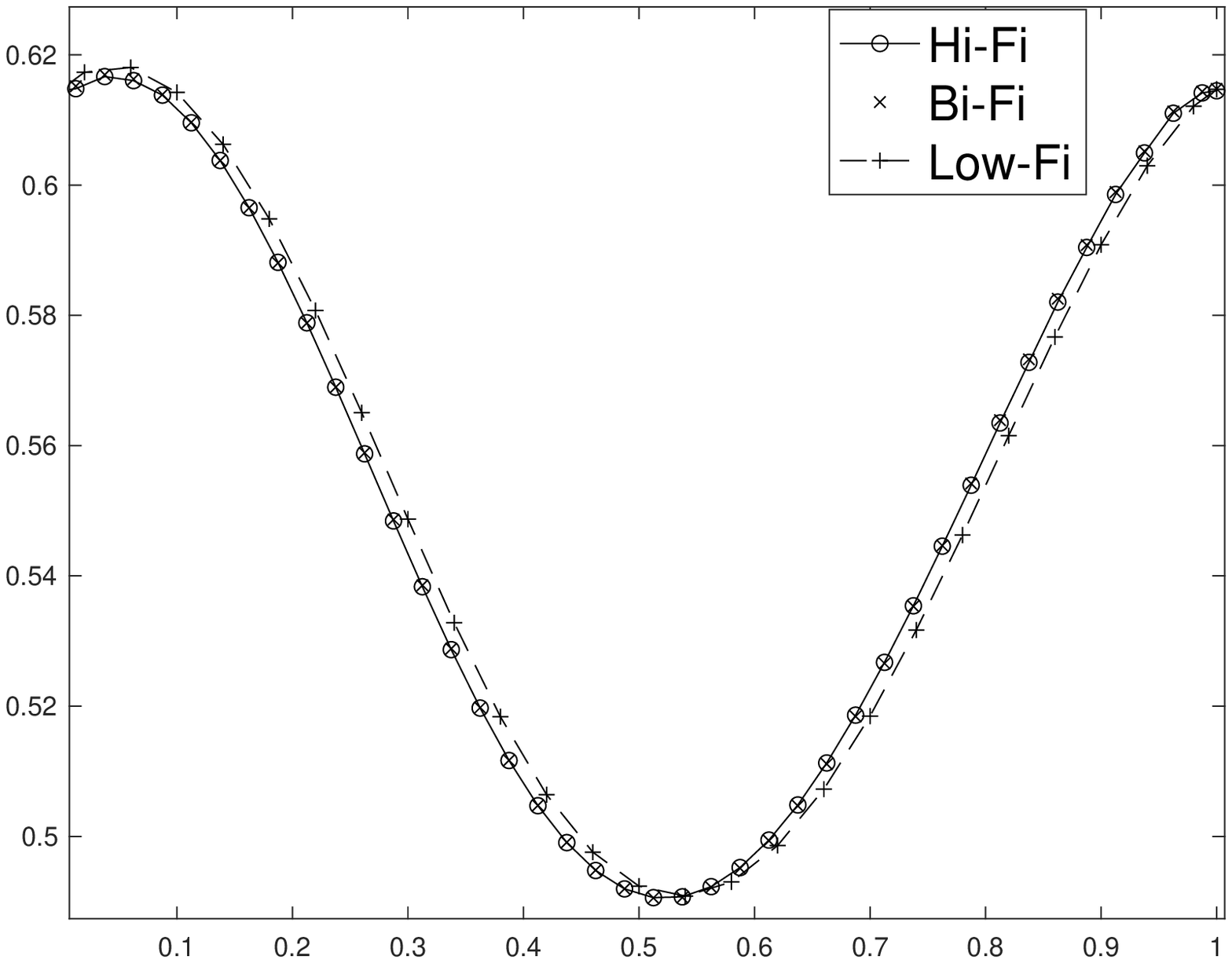}
\includegraphics[scale=0.285]{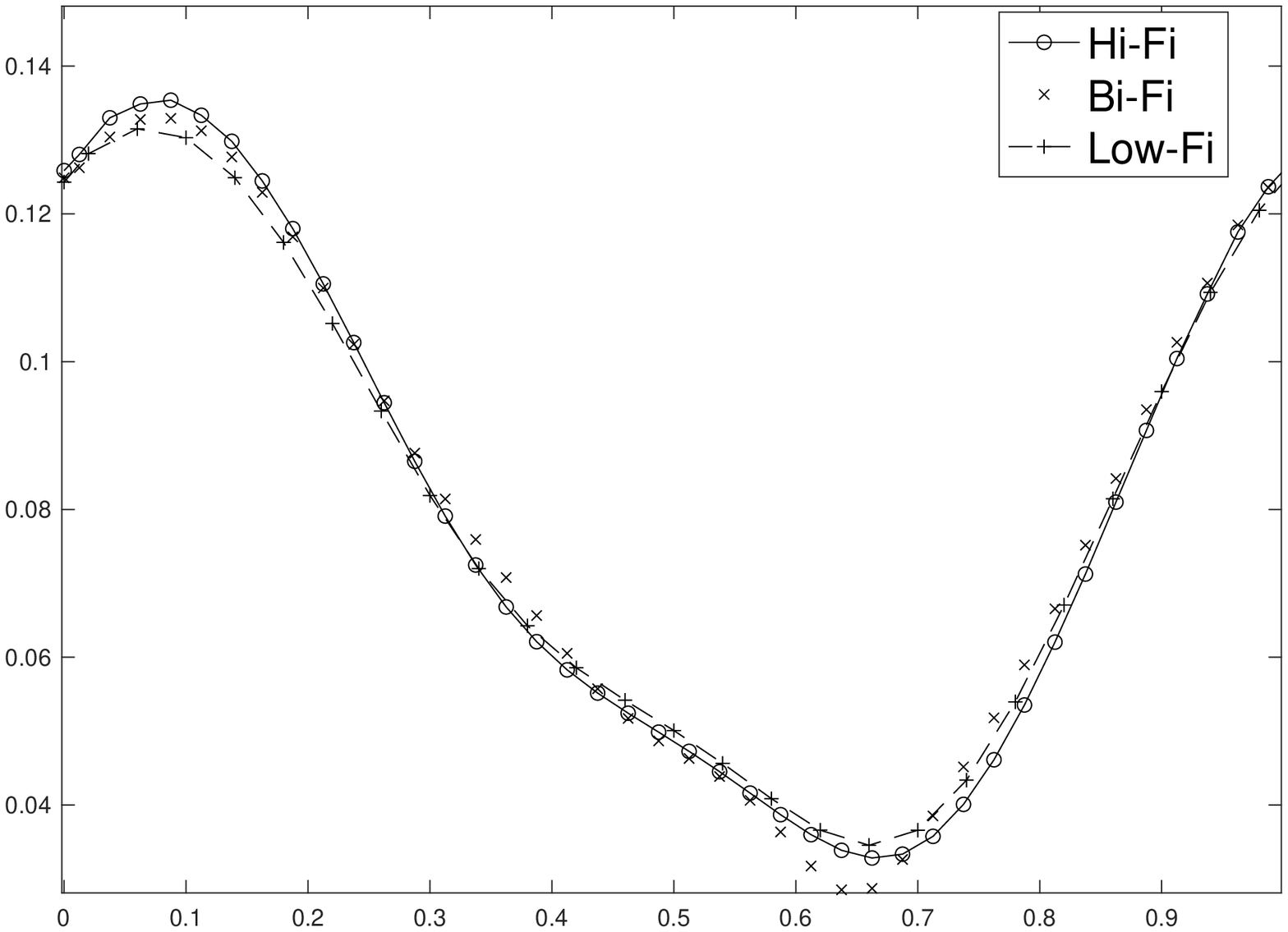}
\includegraphics[scale=0.33]{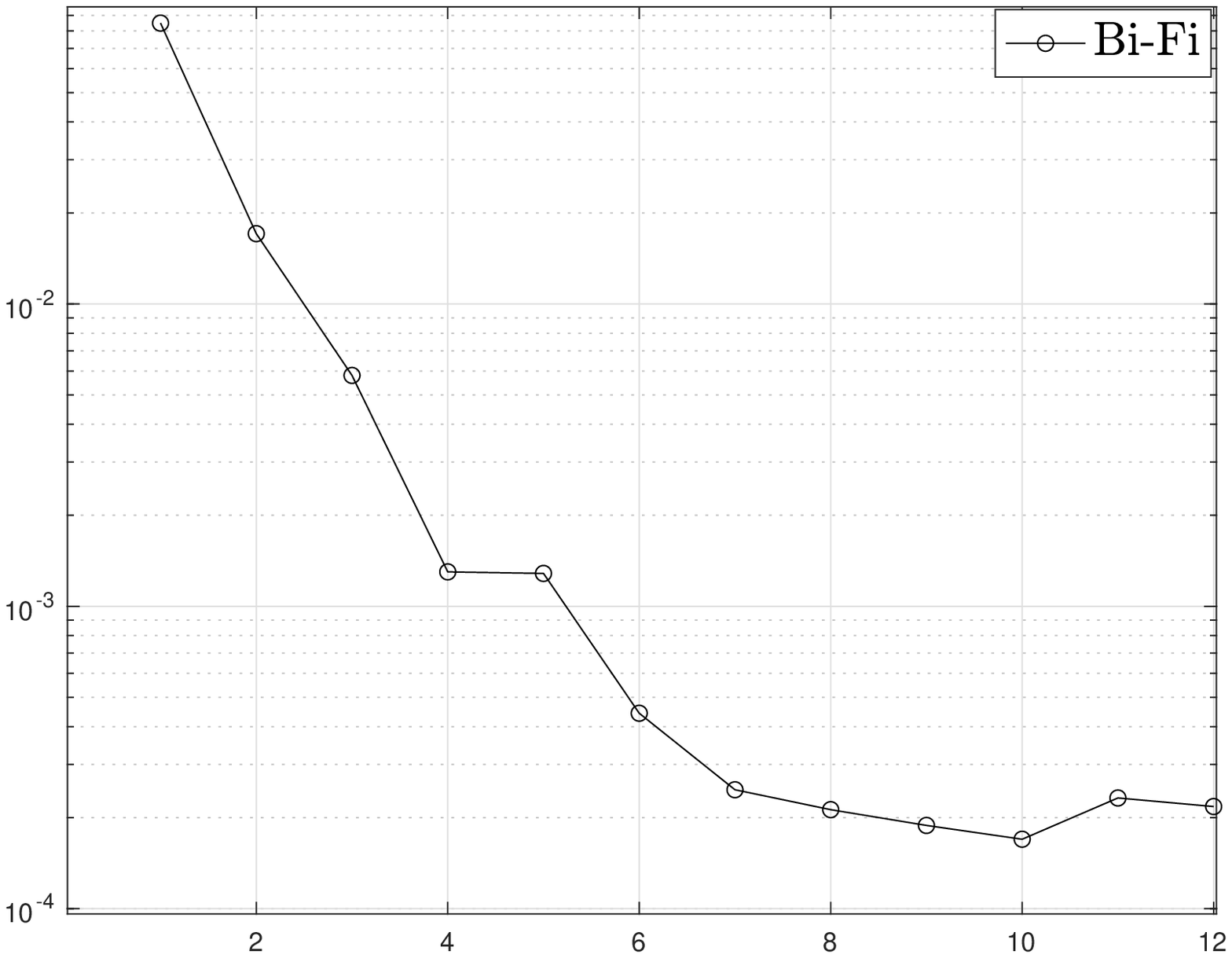}
\includegraphics[scale=0.33]{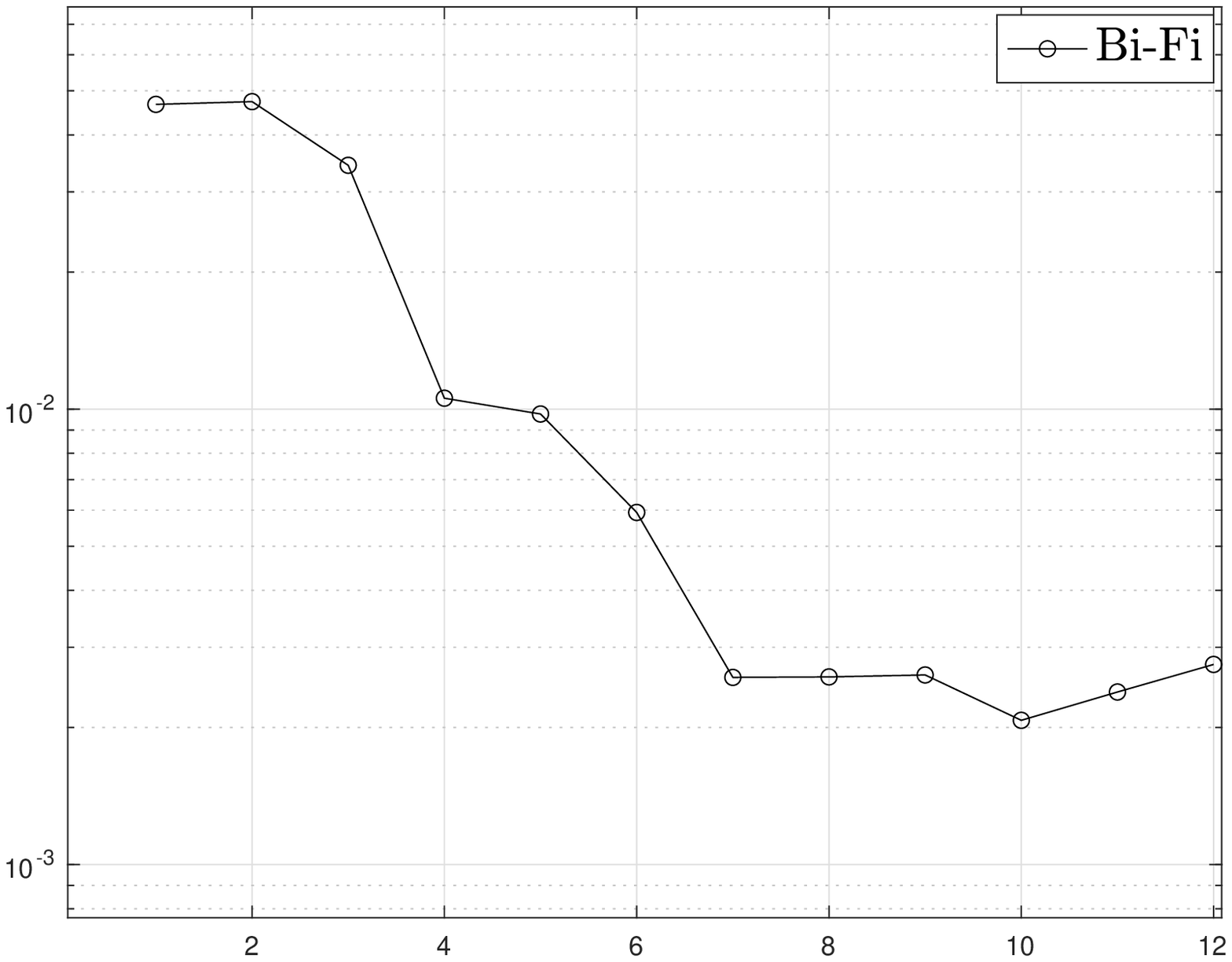}
\caption{Test 2, with $\epsilon=10^{-2}$. The mean (left) and standard deviation (right) of $\overline{r}$, obtained by $r=12$ high-fidelity runs and the sparse grid method with $2243$ quadrature points (crosses, first row). The corresponding errors are also reported (second row).}
\label{Test4_Fig1}
\end{figure}

In this example, we choose $\Delta x=0.025$, $\Delta t=10^{-4}$ in the HF solver, and $\Delta x=0.04$, $\Delta t=2\times 10^{-4}$ in the LF solver.
Numerical solutions and bi-fidelity errors at output time $T=0.02$ with $\epsilon=10^{-2}$ are shown in Fig.~\ref{Test4_Fig1}. 

\begin{figure}[htb]
\centering
\includegraphics[scale=0.3]{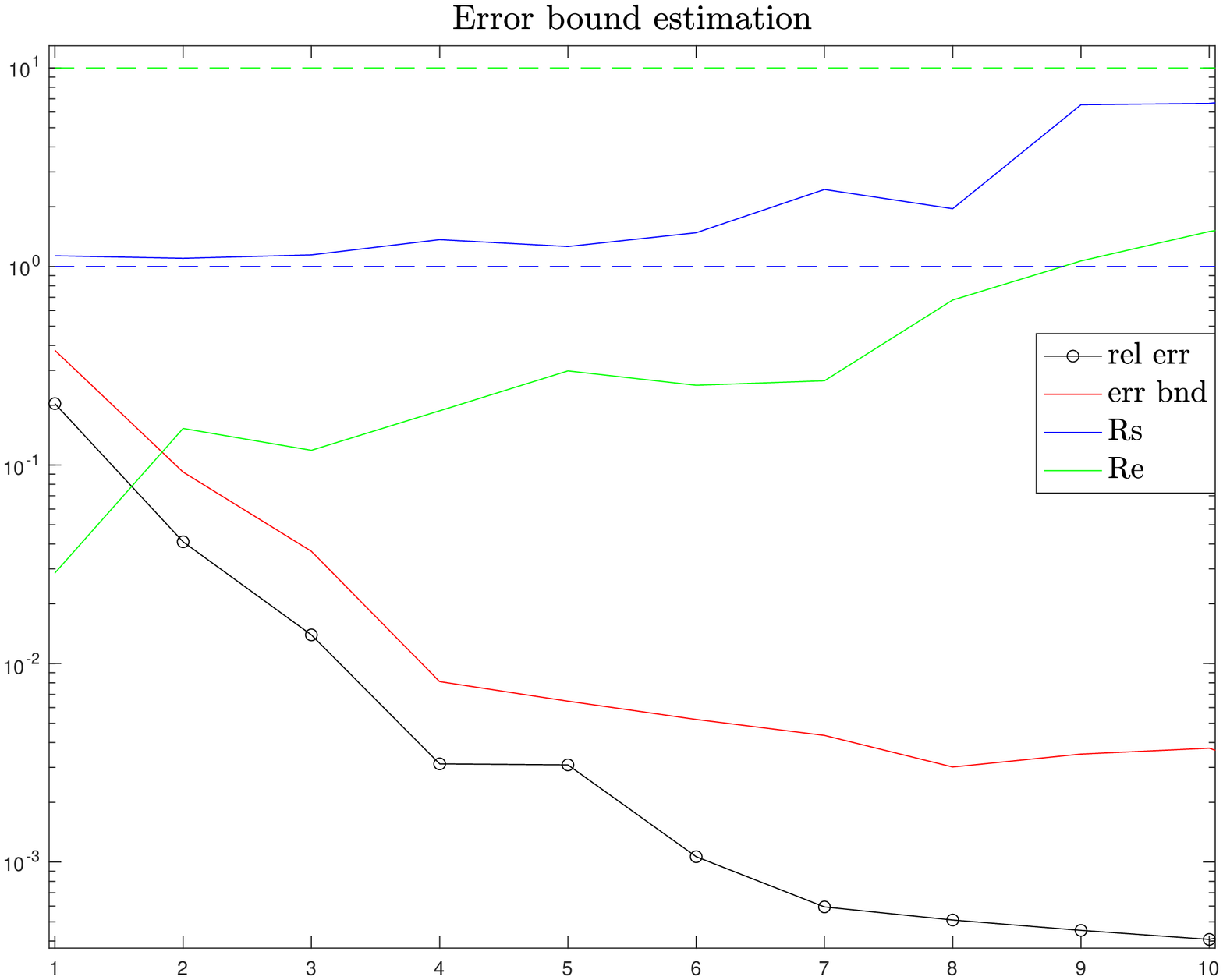}
\caption{Error bound estimation for Test 2.}
\label{Test4_EB}
\end{figure}

From Fig.~\ref{Test4_Fig1} we observe that: while the LF solution misses details of the HF solution near the peaks, the mean and standard deviation of the BF approximations agree well with the HF solution globally, using $r=12$ HF runs. 
A fast exponential decay of errors for the BF approximation is observed from Fig.~\ref{Test4_EB}. With only $12$ HF simulations, the BF mean can reach the accuracy level of less than $\mathcal{O}(10^{-4})$. In Fig.~\ref{Test4_EB}, the practical error bound estimators bound well the true BF errors, with the values of both metrics $R_e$ and $R_s$ are about less than $10$. Similar conclusion as that in Test 1 can be drawn.

\subsection*{Test 3: Riemann problem}
We next consider a Riemann problem assuming again the random coefficient \eqref{Sigma}. Both the boundary condition
\begin{eqnarray*}
	f(x=0,t,v,z)=1+0.4z_1, \quad\mbox{  if }v\ge 0 , \\[2pt]
	 f(x=1,t,v,z)=0, \quad\mbox{  if }v\le 0,
 \end{eqnarray*}
and initial distribution
 \begin{eqnarray*}
 	f(x,t=0,v,z)=1+0.4z_1, \qquad  0\leq x <0.5 , \\[2pt]
 	f(x,t=0,v,z)=0, \qquad  0.5\leq x \leq 1, 
 \end{eqnarray*}
are uncertain. The corresponding B.C. and I.C. for the even and odd parts are 
\begin{equation}
\sigma j=-vr_x .
\end{equation}
thus for $v>0$, 
\begin{equation}
\left. r-\frac{\epsilon}{\sigma} v r_x \right|_{x=0} = 1+0.4z_1, \qquad
\left. r+\frac{\epsilon}{\sigma} v r_x \right|_{x=1} = 0, 
\end{equation}
and
\begin{align*}
 	r(x,t=0,z)&=1+0.4z_1,\quad 0 \leq x <0.5, \\
 	r(x,t=0,z)&=0, \qquad\qquad 0.5\leq x \leq 1, \\
 	j(x,t=0,z)&=0. 
 \end{align*}
In this example, we let the output time $T=0.01$ and $\epsilon=10^{-8}$. 
Set $\Delta t=2\times 10^{-4}$, $\Delta x=0.04$ in the LF solver, and $\Delta t=5\times 10^{-5}$, $\Delta x=0.0125$ in the HF solver. 

\begin{figure}[htb]
\centering
\includegraphics[scale=0.33]{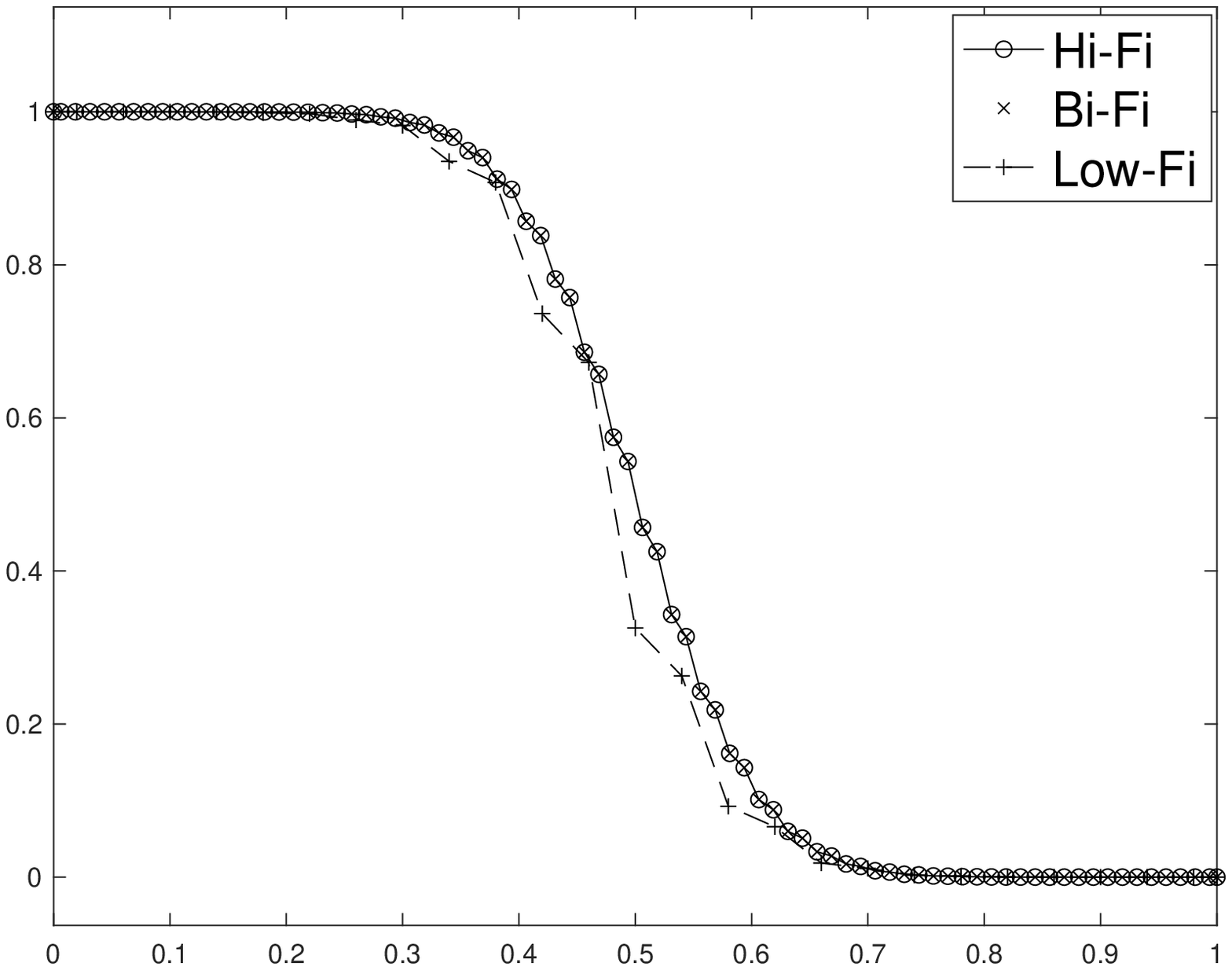} 
\includegraphics[scale=0.33]{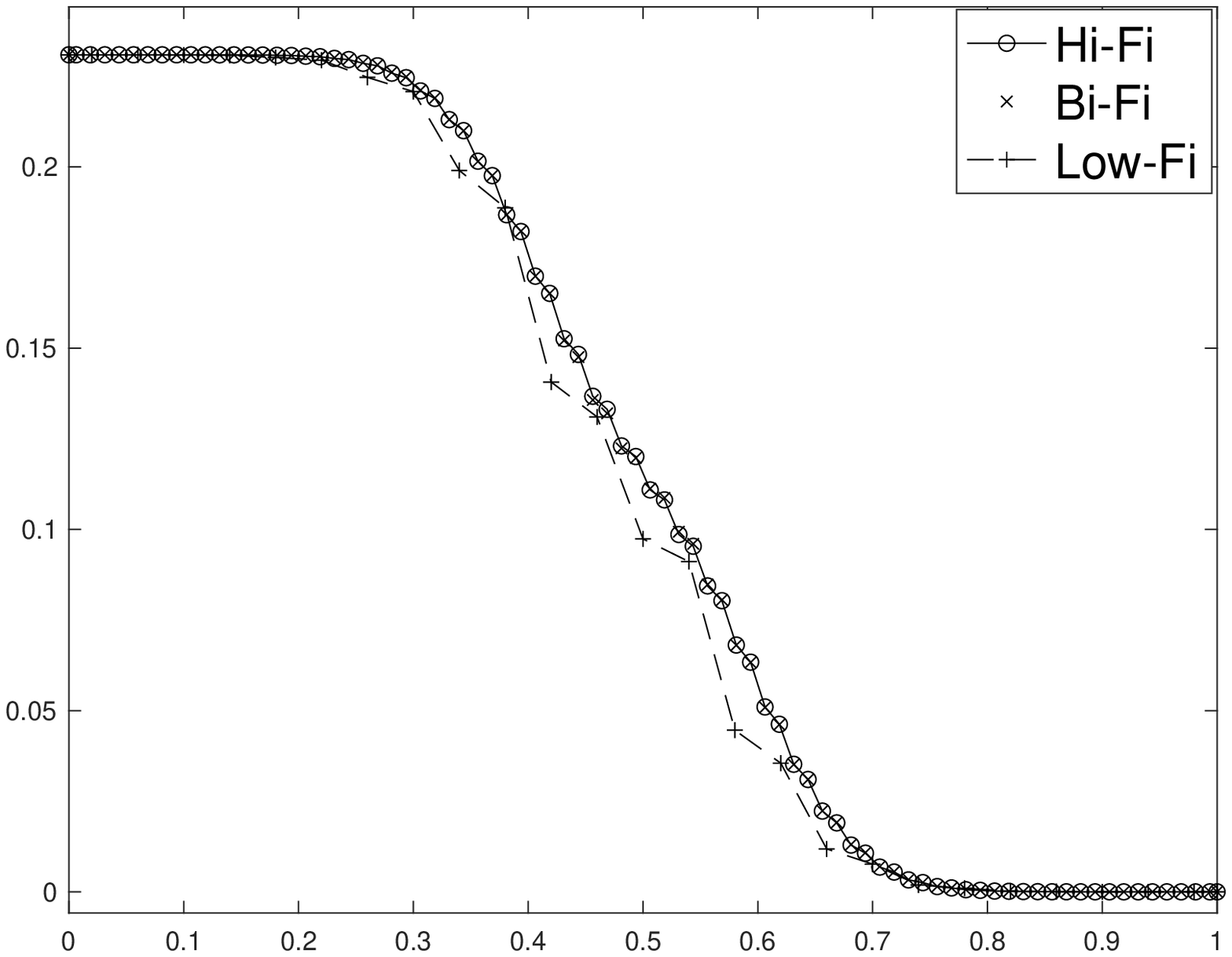}
\caption{Test 3, with $\epsilon=10^{-8}$.
The mean (left) and standard deviation (right) of $\overline{r}$, obtained by $r=12$ high-fidelity runs and the sparse grid method with $2243$ quadrature points (crosses).} 
\label{Test2_Fig1}
\end{figure}

Fig.~\ref{Test2_Fig1} shows the numerical mean and standard deviation of $\bar{r}$ for $\epsilon=10^{-8}$. 
While the LF solutions are not able to capture the detailed information around the transition region, the BF approximation agrees with the HF solutions at a very satisfactory level. 

\begin{figure}[htb]
\centering
\includegraphics[scale=0.33]{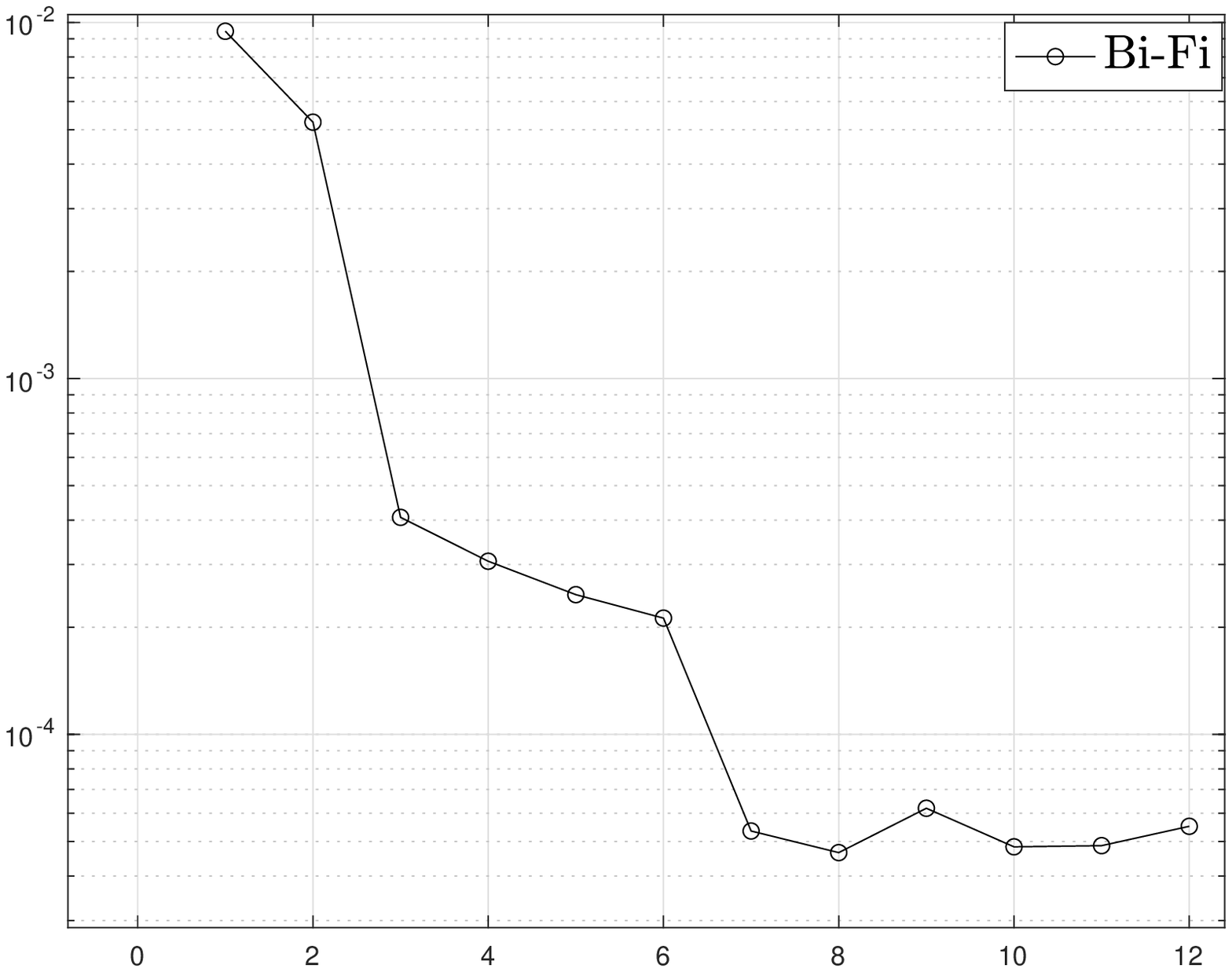}
\includegraphics[scale=0.33]{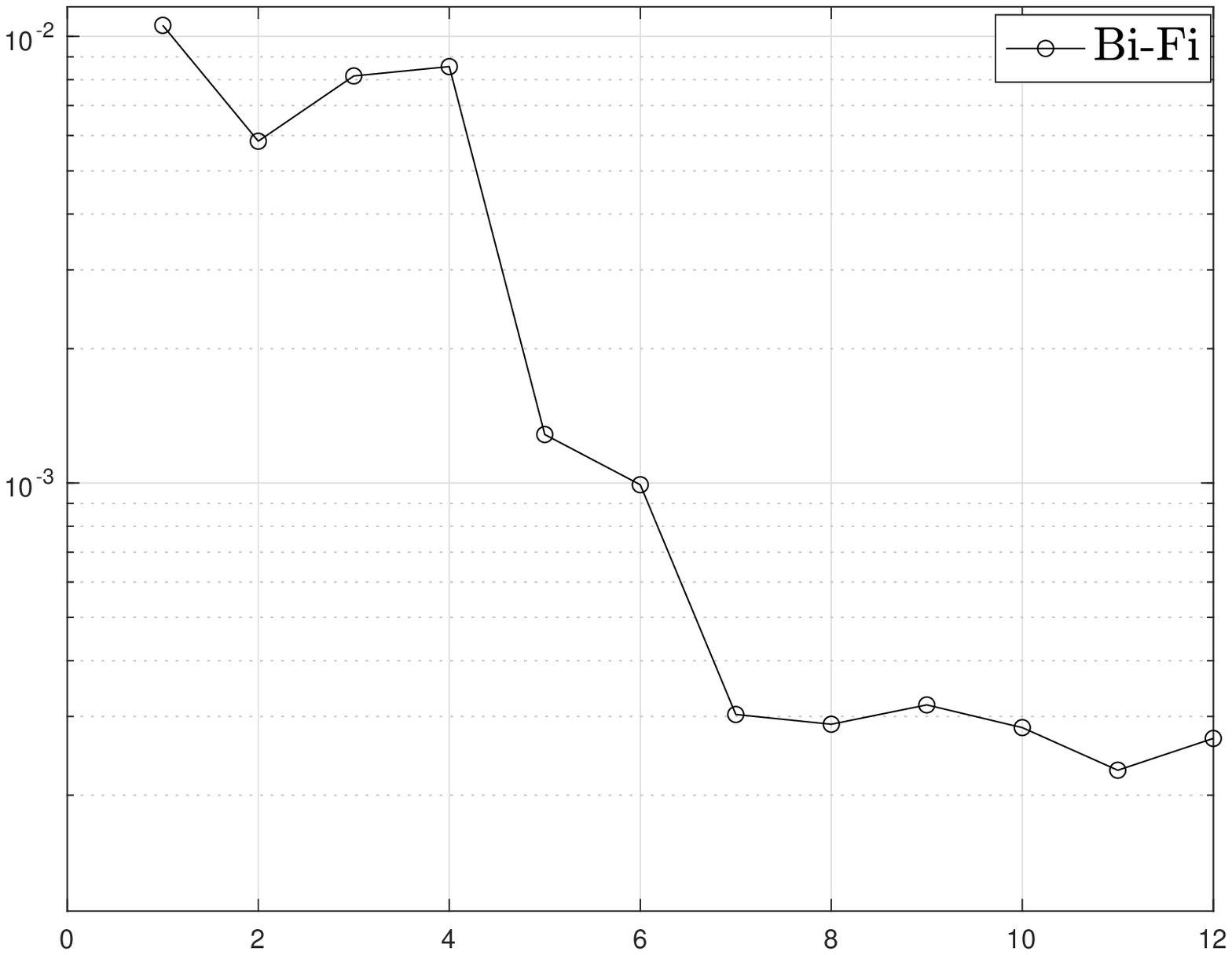}
\caption{Test 3, with $\epsilon=10^{-8}$. 
Errors of the bi-fidelity approximation mean (left) and standard deviation (right) of $\overline{r}$ with respect to the number of high-fidelity runs.}
\label{Test2_Err}
\end{figure}

\begin{figure}[htb]
\centering
\includegraphics[scale=0.4]{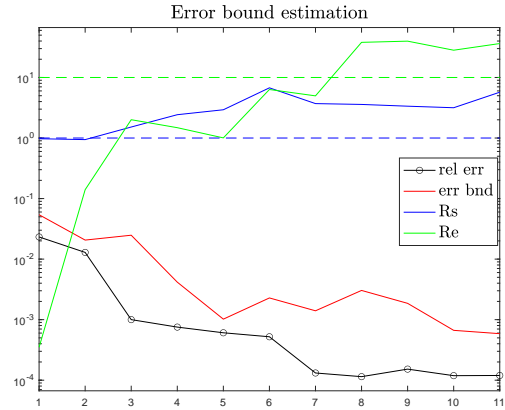}
\caption{Error bound estimation for Test 3.}
\label{Test2_EB}
\end{figure}

From Fig.~\ref{Test2_Err}, it can been seen that the BF errors decay really fast with respect to the selected HF runs. With only $r=12$ HF simulation runs, the BF errors of the mean and standard deviation of $\overline{r}$ reach as small as about $\mathcal{O}(10^{-4})$. Based on Fig.~\ref{Test2_EB}, the practical error bound estimators bound well the true BF errors, with the values of both metrics $R_e$ and $R_s$ around or less than $10$. Again, further increasing the HF samples after $r=8$ does not help in improving the quality of BF approximation, which is consistent with our expectations. 

\subsection*{Test 4: Mixed regime test}
\begin{figure}[h!]
\centering
\includegraphics[scale=0.4]{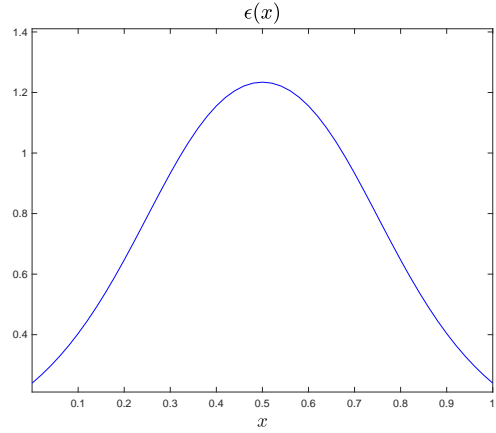}
\caption{Illustration of the spatially varying $\epsilon(x)$ over the physical domain for Test 4.}
\label{mixedeps}
\end{figure}

\begin{figure}[htb]
\centering
\includegraphics[scale=0.33]{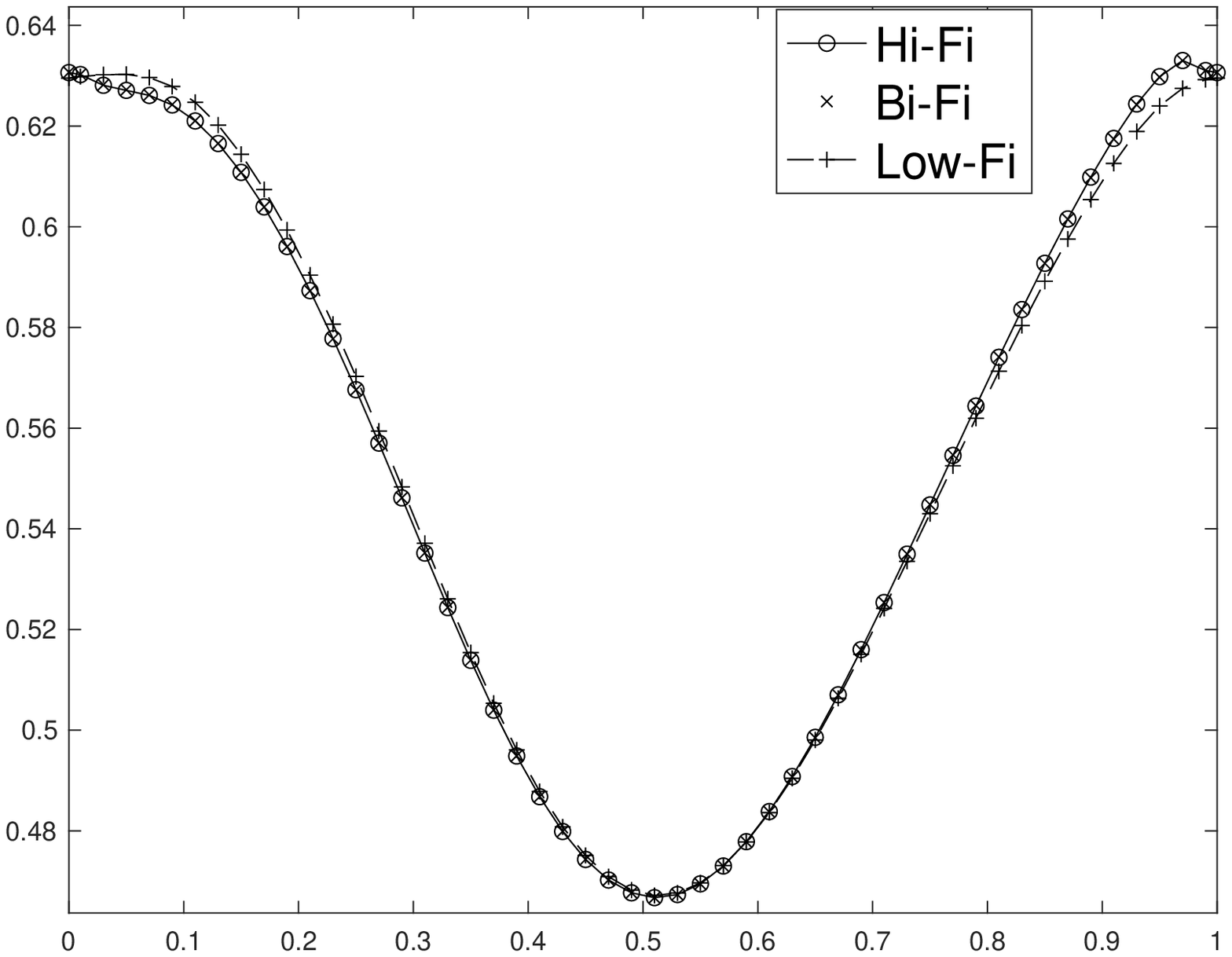}
\includegraphics[scale=0.33]{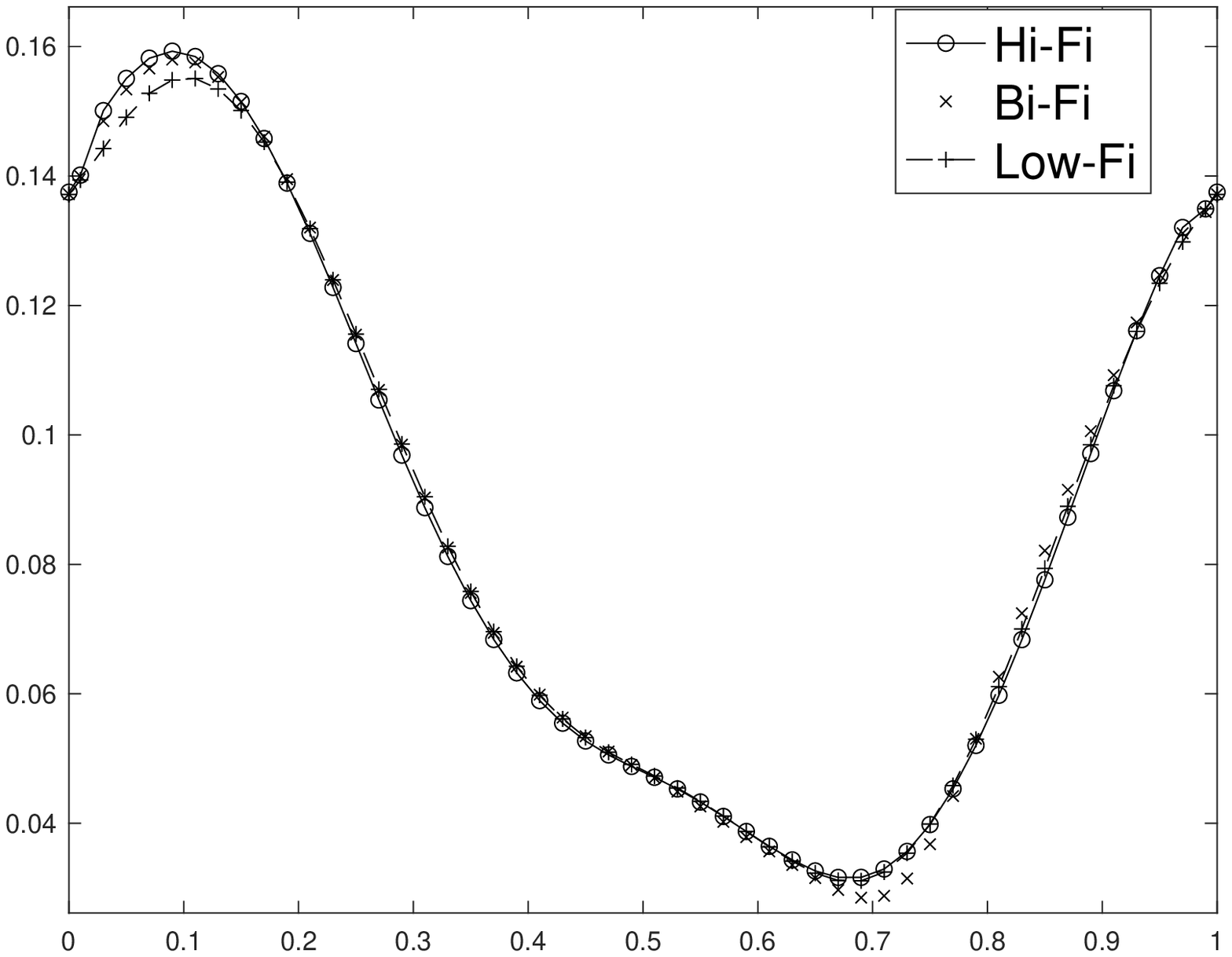}
\caption{Test 4. The mean (left) and standard deviation (right) of $\overline{r}$, obtained by $r=12$ high-fidelity runs and the sparse grid method with $2243$ quadrature points (crosses).}
\label{Test3_Fig1}
\end{figure}

\begin{figure}[htb]
\centering
\includegraphics[scale=0.33]{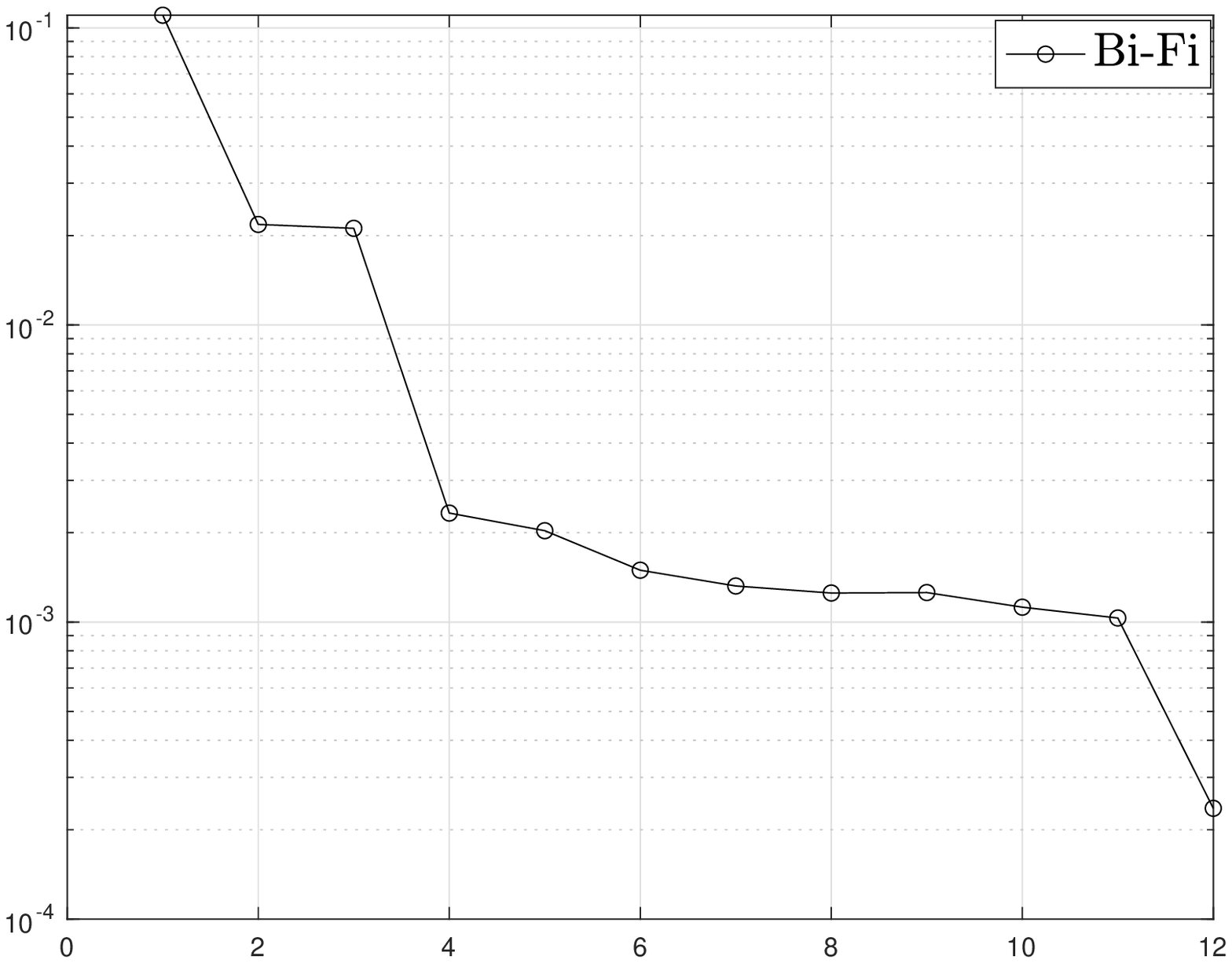}
\includegraphics[scale=0.33]{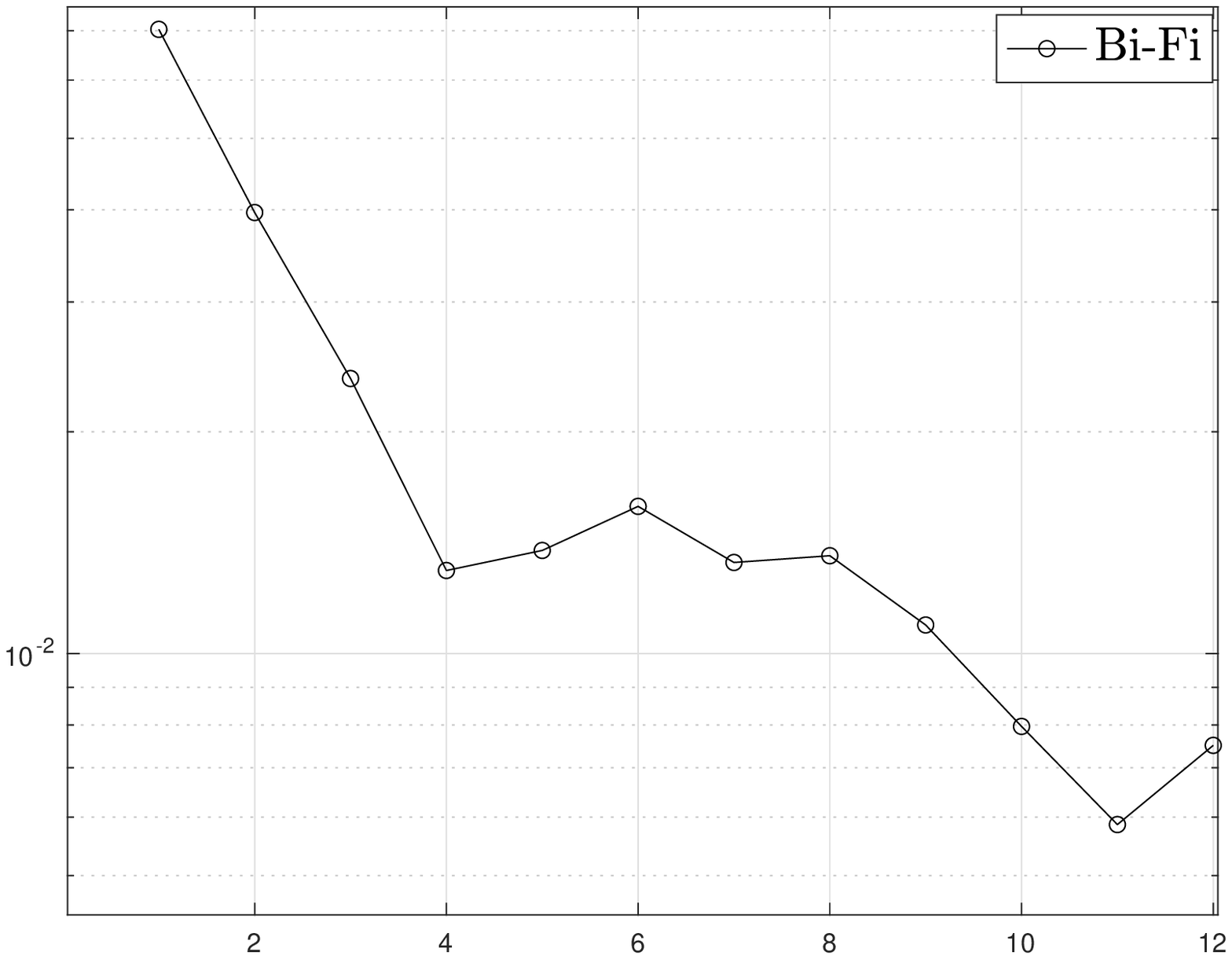}
\caption{Test 4. Errors of the bi-fidelity approximation mean (left) and standard deviation (right) of $\overline{r}$ with respect to the number of high-fidelity runs.}
\label{Test3_Err} 
\end{figure}

\begin{figure}[htb]
\centering
\includegraphics[scale=0.4]{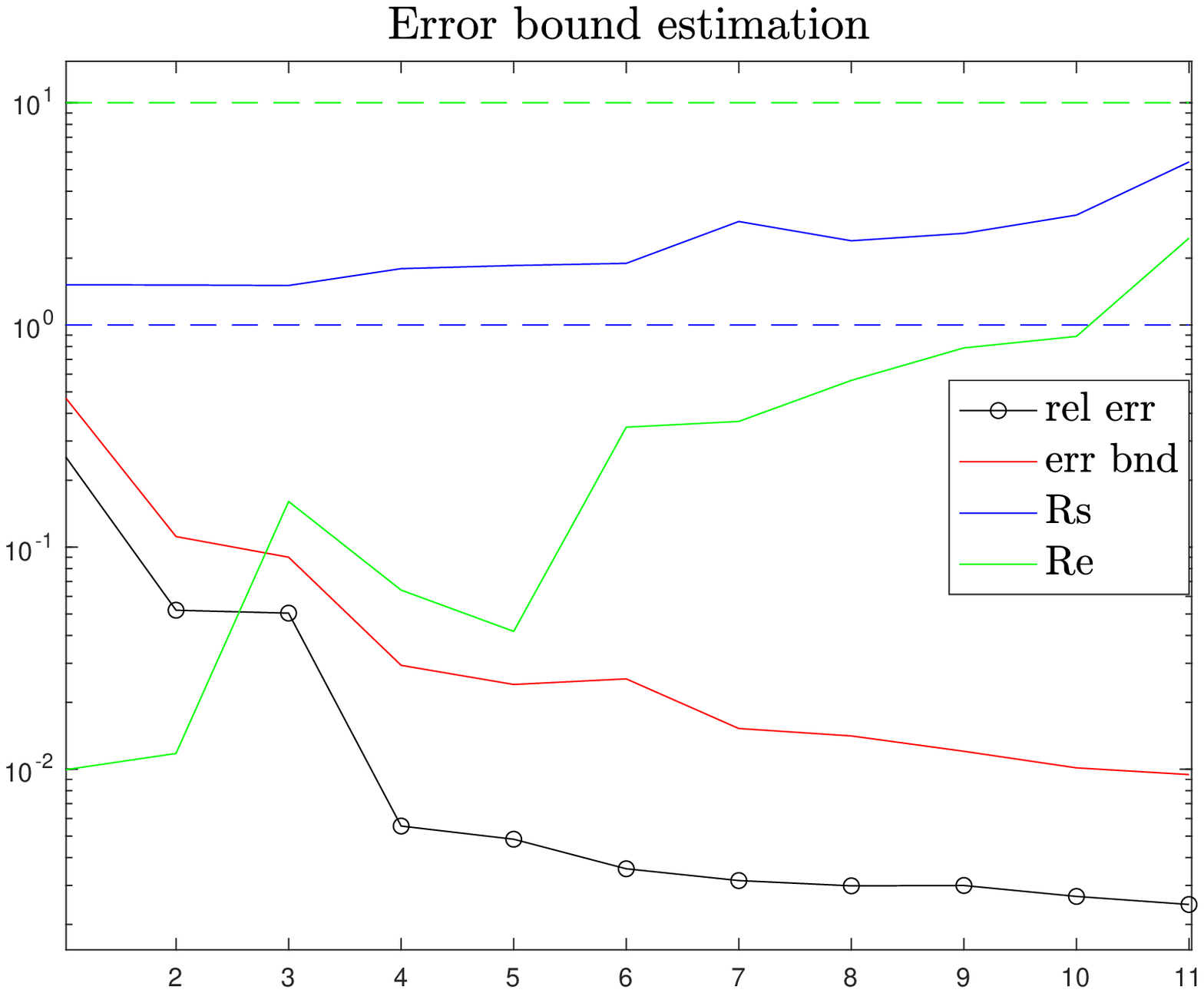}
\caption{Error bound estimation for Test 4.}
\label{Test3_EB} 
\end{figure}

In Test 4, we investigate the performance of the BF approximation in a mixed regime, which is a benchmark test for multiscale kinetic problems. 
In this case, the Knudsen number is spatially dependent and is given by
\begin{equation}
\epsilon^2(x)= 10^{-8} + \left[\tanh\left(1-\frac{11}{2}(x-0.5)\right)+\tanh\left(1+\frac{11}{2}(x-0.5)\right)\right]. 
\end{equation}
As shown in Fig.~\ref{mixedeps}, the values of $\epsilon(x)$ vary smoothly from $\mathcal O(10^{-1})$ to $O(1)$.   
Assume a deterministic cross-section $\sigma=1$, {and the same initial data as \eqref{Test2_IC}.}
Periodic boundary conditions are considered here. We use $\Delta t=10^{-4}$ in the LF, $\Delta t=5\times 10^{-5}$ in the HF solver, and $\Delta x=0.02$ in both solvers.
Numerical results and error bound estimation at time $T=0.01$ are shown in Fig.~\ref{Test3_Fig1} -- \ref{Test3_EB}.

In contrast to the poor approximation of LF solutions in some regions, 
the mean and standard deviation of the BF approximations match very well with the HF solution across the whole range of $\epsilon$, {\color{black}suggesting the GT model can capture the important variations of the linear transport model in the uncertainty parameter space across multiple scales}. Again, a fast exponential decay of the errors for the BF approximation is observed. With only $r=12$ HF simulations, the BF 
mean can reach the accuracy of about $\mathcal{O}(10^{-4})$, which is quite remarkable compared to other non-intrusive sampling methods such as the Monte-Carlo based or sparse-grid type methods. 

In Fig.~\ref{Test3_EB}, both two metrics $R_s$ and $R_e$ are less than $10$. 
It is worthy to mention that the true error curve (black line with dots)
are computed based on the HF solutions of the entire test sets, while the error bound estimations (red line) only depends on the existing pre-selected training HF data, suggesting the effectiveness of the error bound estimation. 
In practical applications, the suggested evaluation metrics and error bound estimation could indeed provide an admissible and easy way to evaluate the performance of the BF method a priori. 


\subsection*{Test 5: A discontinuous cross-section test}

In the last test case, we consider the cross-section to be uncertain and discontinuous: 
\begin{equation}
\sigma(x,z)=\left\{
\begin{split}
& 1 + 4 \sum_{i=1}^d \frac{1}{(i\pi)^2}\cos{(2\pi i x)} z_i, \qquad 0 \leq x \leq 0.5, \\[4pt]
& 0.2, \qquad 0.5<x<1, 
\end{split}
\right.
\end{equation}
where $d=5$. The initial distribution is given by 
\begin{equation} f(t=0,x,v) = 
\displaystyle \frac{1}{2\pi\xi}\exp\left(-\frac{(x-0.5)^2}{2\xi}\right), \qquad 
\xi=0.01. 
\end{equation}
Periodic boundary condition is assumed. We set $\epsilon^2=0.001$ and the output time $T=0.01$. 
We use $\Delta x=0.025$, $\Delta t=2\times 10^{-4}$ in the LF solver, and  $\Delta x=0.02$, $\Delta t=2/3\times 10^{-4}$ in the HF solver. 

\begin{figure}[htb]
\centering
\includegraphics[scale=0.33]{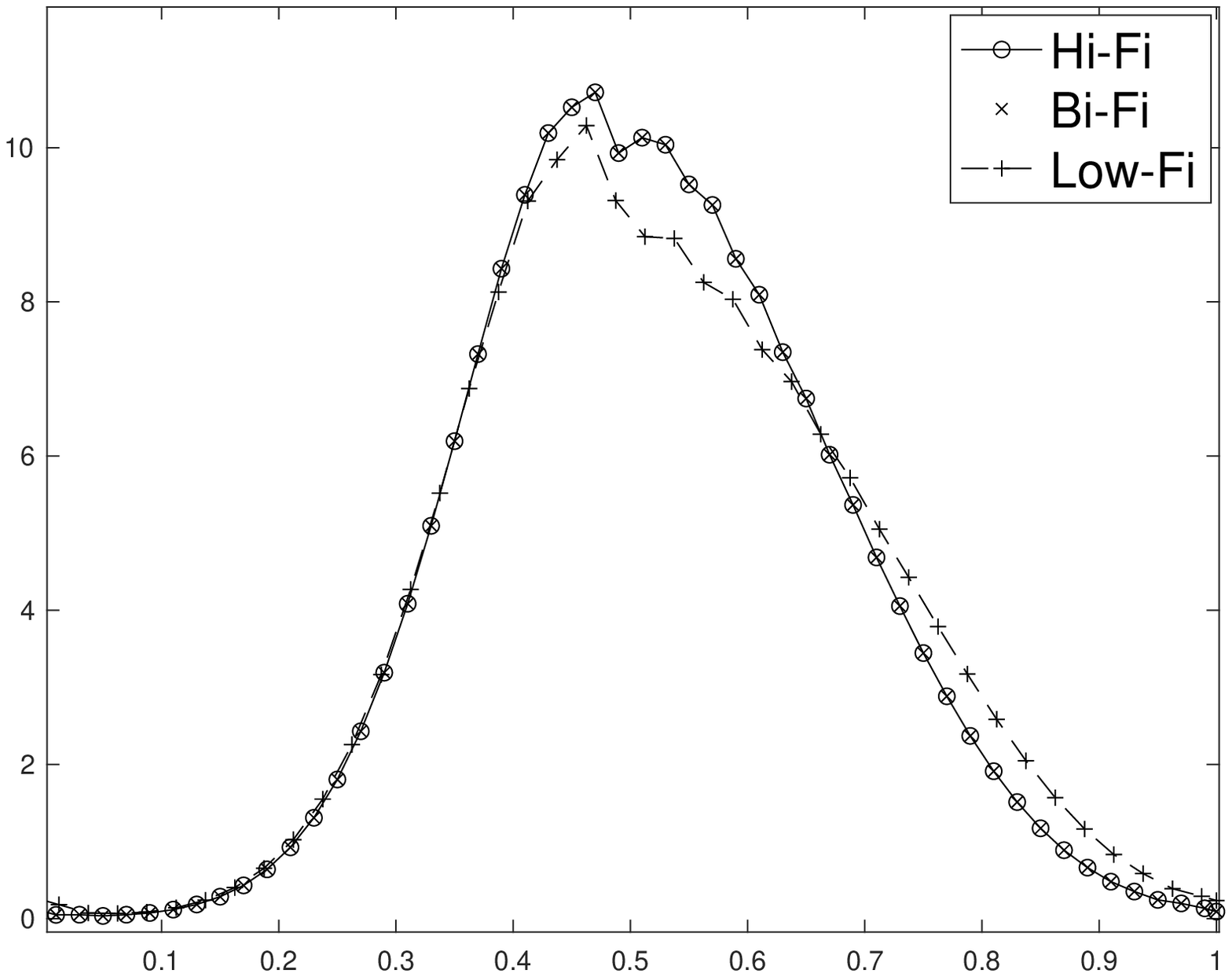}
\includegraphics[scale=0.33]{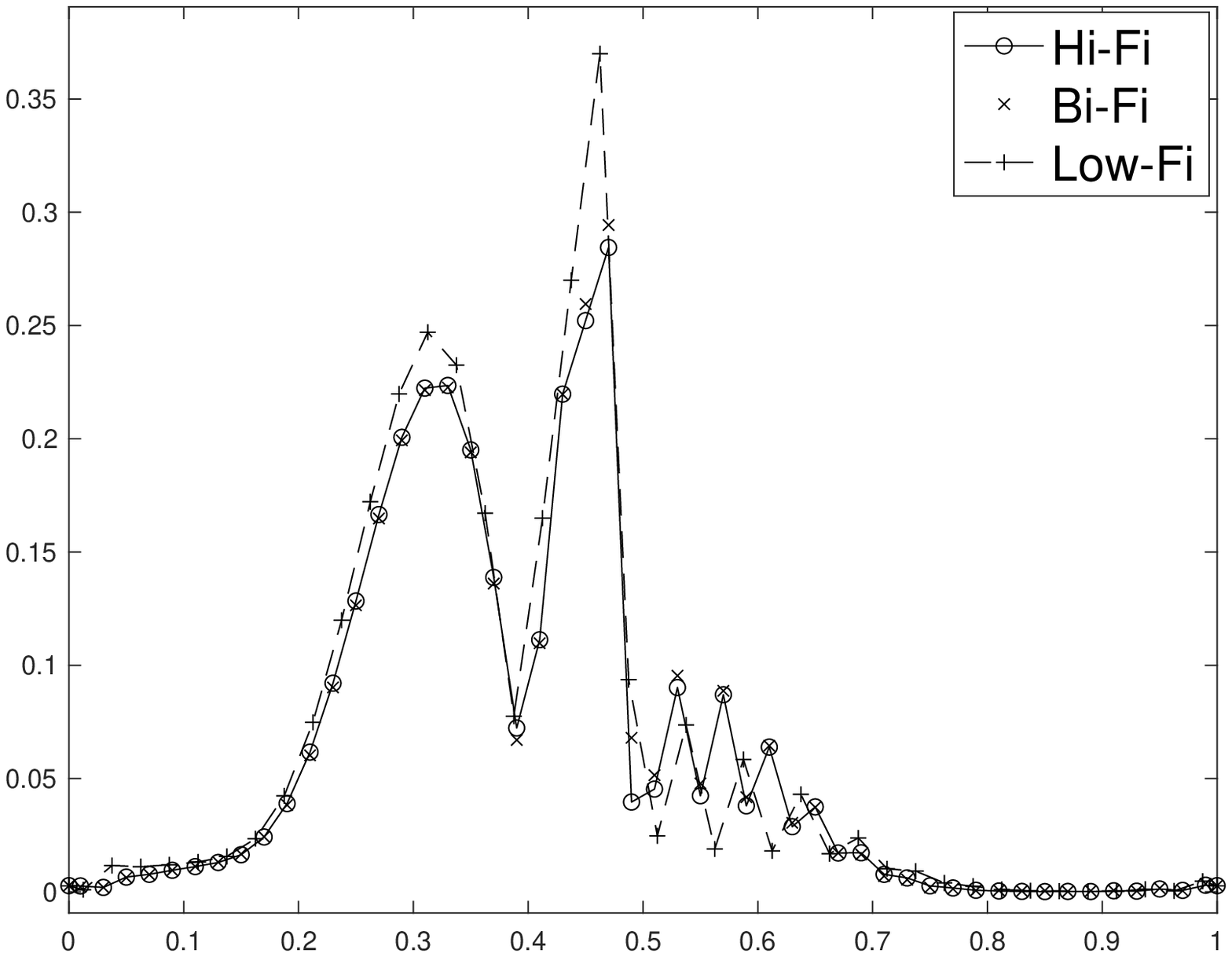}
\includegraphics[scale=0.33]{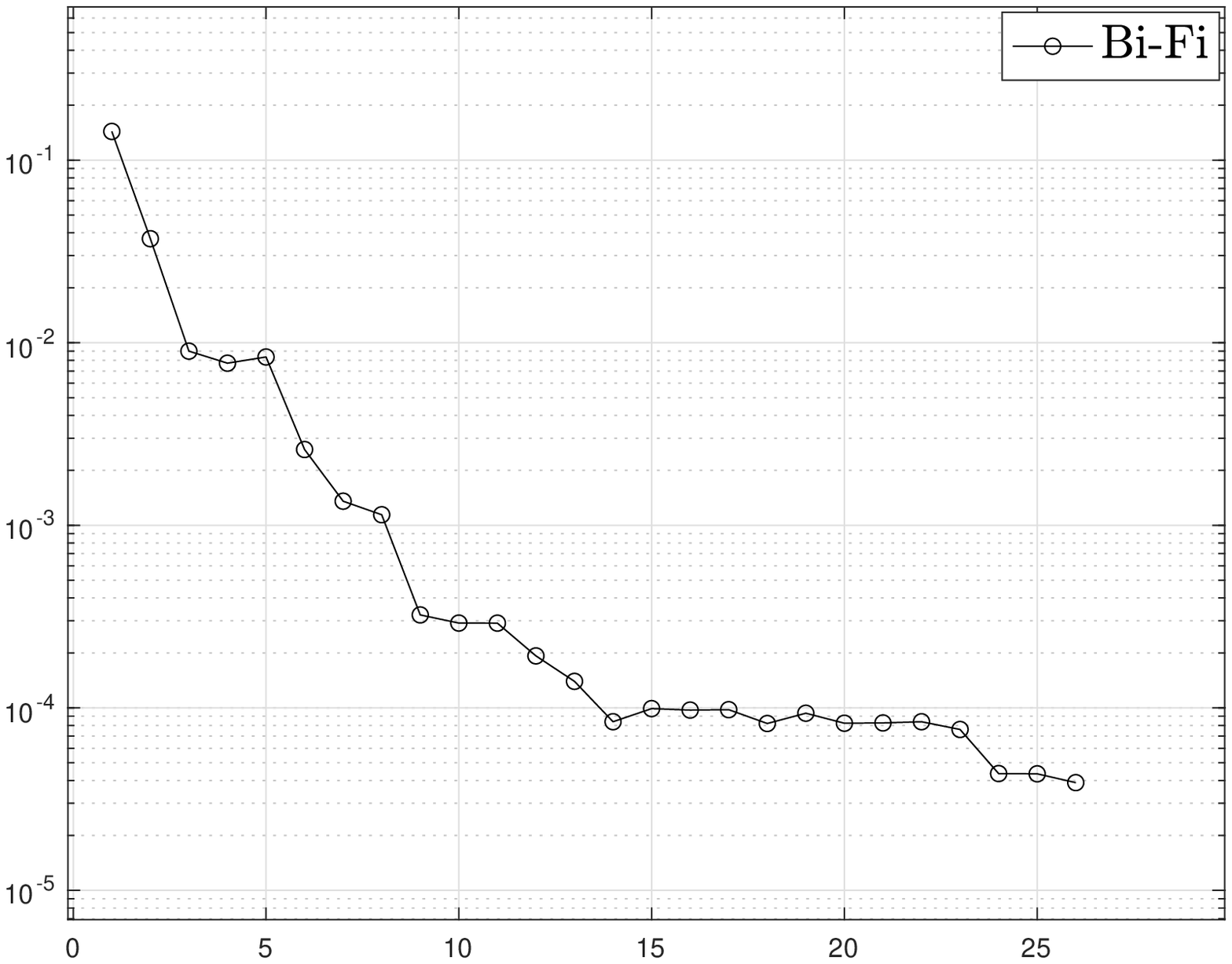}
\includegraphics[scale=0.33]{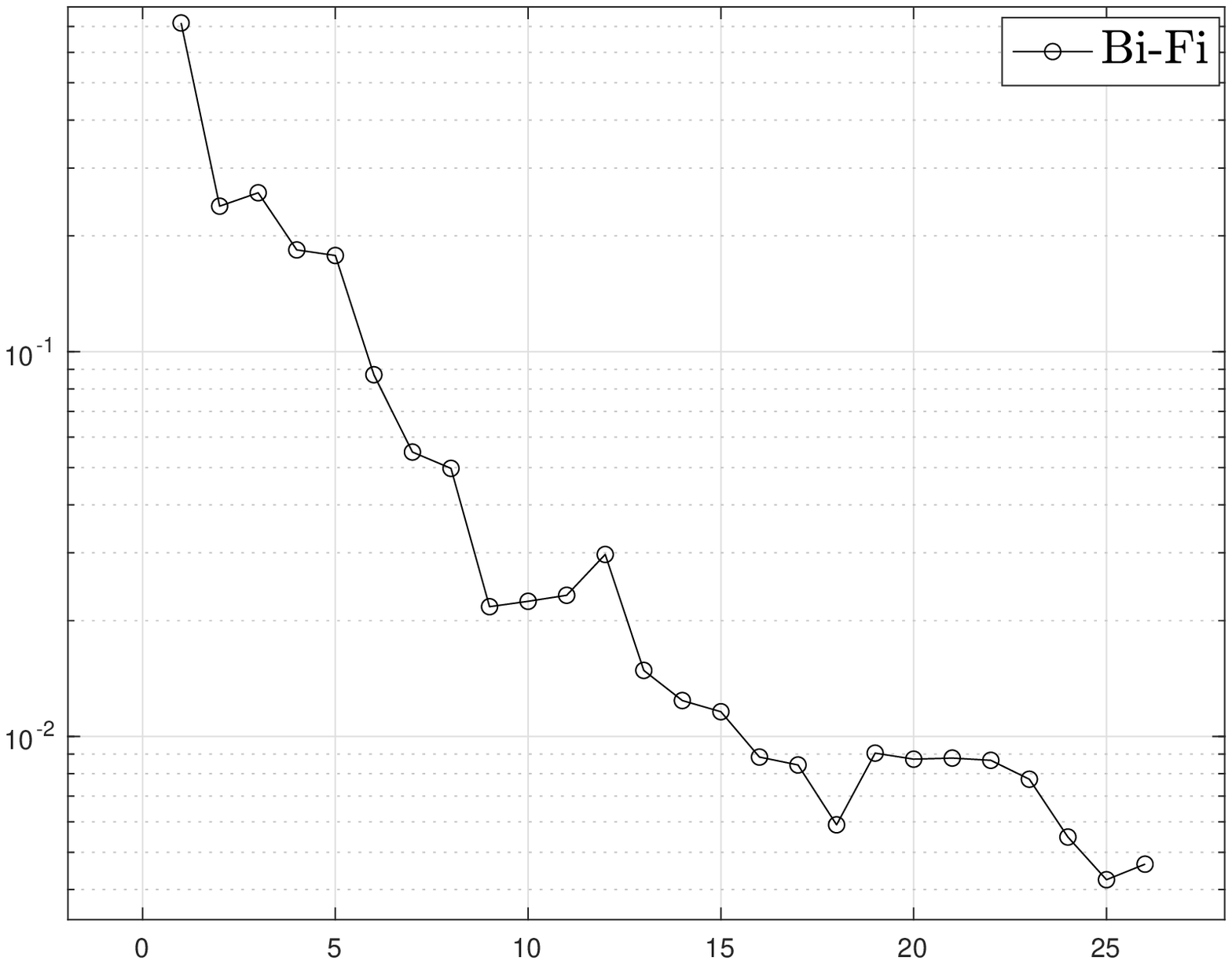}
\caption{Test 5. The mean (left) and standard deviation (right) of $\overline{r}$, obtained by $r=25$ high-fidelity runs and the sparse grid method with $2243$ quadrature points (crosses, first row). The corresponding errors are also reported (second row).}
\label{Test5}
\end{figure}

\begin{figure}[htb]

\centering
\includegraphics[scale=0.4]{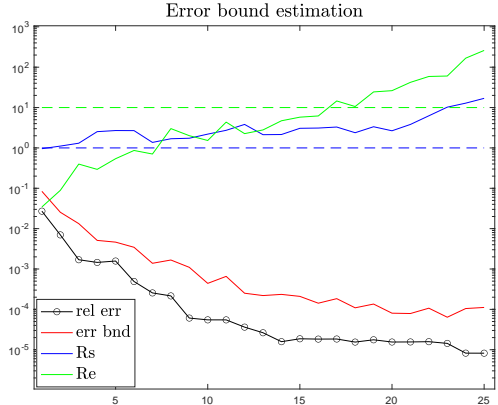}
\caption{Error bound estimation for Test 5.}
\label{Test5_EB} 
\end{figure}

From Fig.~\ref{Test5}, 
we observe that the mean and standard deviation of the BF approximations match quite well with the HF solutions, particularly around the discontinuity. The BF approximations again enjoy a fast exponential decay of the errors. With only $15$ HF simulations, the BF mean reaches the accuracy of about $\mathcal O(10^{-4})$. In Fig.~\ref{Test5_EB}, the practical error bound estimators bound the true BF errors well, while the values of both metrics $R_e$ and $R_s$ are less than $10$. 
After $r=15$, if one increases the HF samples, it will not help to improve the quality of BF approximations. 

\section{Conclusion}

In this paper, we investigate the applicability and performance of the bi-fidelity stochastic collocation method to quantify the uncertainties of multiscale linear transport equations with multi-dimensional random parameters. The sources of uncertainties considered are multidimensional and include the collision kernel or initial and boundary conditions. The Goldstein-Taylor model, which is 
a simpler model with discrete velocities compared to the linear transport model under study, has been chosen as the low-fidelity model. Both models share the same diffusion limit. 
Our numerical examples demonstrate that the proposed bi-fidelity method works effectively across a large range of regimes in diffusive, kinetic and mixed regimes. This suggest that the present approach is capable to reach a uniform accuracy with respect to the multiple scale involved.
Furthermore, an empirical error bound estimation is computed to access the quality of the bi-fidelity approximation. Further researches will consider the extension of the present approach to transport equations occurring in biomathematics, like chemotaxis and epidemiology \cite{JH1, boscheri2020, Ber, Ber2}. 

\section*{Acknowledgments}
X. Zhu was supported by the Simons Foundation (504054). L. Liu was supported by the start-up funding from CUHK and Early Career Scheme  by Research Grants Council of Hong Kong (24301021). L. Pareschi was partially supported by MIUR (Ministero dell’Istruzione, dell’Università e della
Ricerca) PRIN 2017, project “\textit{Innovative numerical methods for evolutionary partial differential equations and applications}”, code 2017KKJP4X and by GNCS (Gruppo Nazionale per il Calcolo Scientifico) of INdAM (Istituto Nazionale di Alta Matematica).

\end{document}